\documentclass{article}

\usepackage{amsmath,amssymb,bm}
\usepackage{cite,url,color}

\usepackage[top=2.5cm,bottom=2.5cm,left=2cm,right=2cm]{geometry}

\newtheorem{definition}{Definition}
\newtheorem{theorem}{Theorem}
\newtheorem{lemma}{Lemma}
\newtheorem{corollary}{Corollary}
\newtheorem{algorithm}{Algorithm}
\newtheorem{example}{Example}
\newtheorem{remark}{Remark}
\newenvironment{proof}[1][Proof]{\par\noindent\textit{#1}: }{\hfill$\blacksquare$\vskip 0.5\baselineskip}
\newenvironment{solution}[1][Solution]{\par\noindent\textit{#1}: }{\hfill$\square$\vskip 0.5\baselineskip}
\newenvironment{explain}[1][Explain]{\par\noindent\textit{#1}: }{\hfill$\square$\vskip 0.5\baselineskip}

\usepackage[normalem]{ulem}

\usepackage{arydshln}

\begin{document}

\title{Null Polynomials modulo $m$}
\author{Shujun Li\\\url{http://www.hooklee.com}}
\date{\today}

\maketitle

\begin{abstract}
This paper studies so-called ``null polynomials modulo $m$", i.e.,
polynomials with integer coefficients\footnote{In this paper, we
always call them ``integer polynomial" in short. Note that sometime
another name ``integral polynomial" is used \cite[Sec.
7.2]{HW:NumberTheory1979}. However, we prefer to ``integer
polynomial" to avoid confusion with the word ``integral" as an
adjective (see \cite{IntegerPolynomial}).} that satisfy $f(x)\equiv
0\pmod m$ for any integer $x$. The study on null polynomials is
helpful to reduce congruences of higher degrees modulo $m$ and to
enumerate equivalent polynomial functions modulo $m$, i.e.,
functions over $\mathbb{Z}_m=\{0,\cdots,m-1\}$ generated by integer
polynomials. The most well-known null polynomial is $f(x)=x^p-x$
modulo a prime $p$.

After pointing out that null polynomials modulo a composite can be
studied by handling null polynomials modulo each prime power, this
paper mainly focuses on null polynomials modulo $p^d$ ($d\geq 1$). A
typical monic null polynomial of the least degree modulo $p^d$ is
given for any value of $d\geq 1$, from which one can further
enumerate all null polynomials modulo $p^d$. The most useful result
obtained in this paper are Theorem
\ref{theorem:H-is-least-degree-monic-null-poly} in Sec.
\ref{subsection:null-poly-pd} and its derivative -- Theorem
\ref{theorem:null-poly-enumerate-all} in Sec.
\ref{subsection:null-poly-pd-enumerate-all}. The results given in
Sec. \ref{subsection:null-poly-2-I3} form a basis of the induction
proofs given in Sec. \ref{subsection:null-poly-pd}. However, if you
do not care how the proofs in Sec. \ref{subsection:null-poly-pd}
were established, you can simply skip Sec.
\ref{subsection:null-poly-2-I3}. Theorems
\ref{theorem:base-null-poly-feature} and
\ref{theorem:base-null-poly-hierarchy} are very important for the
proof of Theorem \ref{theorem:H-is-least-degree-monic-null-poly} and
should be paid more attention.

\textit{Note: After finishing this draft, we noticed that some
results given in this paper have been covered in Kempner's papers
\cite{Kempner:PolyResidue:TAMS1921a, Kempner:PolyResidue:TAMS1921b}.
Since we use a different way to obtain the results, this work can be
considered as an independent and different proof. For a brief
introduction to Kempner's proof, see the Appendix of this paper.}
\end{abstract}

\tableofcontents

\section{Introduction}

In this paper, we study integer polynomials that satisfy $f(x)\equiv
0\pmod m$ for any integer $x$. We call such polynomials ``null
polynomials modulo $m$", due to the fact that they generate nothing
meaningful modulo $p^d$. Till now, I have not found an existing name
of such polynomials. If you know one, please let me know and
recommend a paper or book for reference. Thanks in advance for your
help.

When $p$ is a prime, there is a most well-known null polynomial
modulo $p$: $f(x)=x^p-x$. However, actually it is not the simplest
one, since one should use the Fermat's little theorem to derive that
it is a null polynomial modulo $p$. In my opinion, $\forall
m\in\mathbb{Z}$, the most simplest null polynomial modulo $m$ is
$f(x)=\prod_{i=0}^{m-1}(x-i)$. When $p$ is a prime, we have
$\prod_{i=0}^{p-1}(x-i)\equiv x^p-x\pmod p$ \cite[Theorem
112]{HW:NumberTheory1979}, where ``$\equiv$" means that the
coefficients of the two polynomials are congruent modulo $p$.
However, when $m$ is not a prime, $f(x)=\prod_{i=0}^{m-1}(x-i)$ is
generally not a null polynomial of the least degree modulo $m$. For
example, when $m=p^d$ ($d\geq 2$), it is obvious that
$f(x)=p^{d-1}(x^p-x)$ is a null polynomial of degree $p$ modulo
$p^d$ and $f(x)=(x^p-x)^d$ is a \textbf{monic} null polynomial of
degree $pd$ modulo $p^d$. Then, we have a question: can we find the
least degree of all (monic) null polynomials modulo any integer $m$?
This paper gives an affirmative answer to this question (Theorems
\ref{theorem:H-is-least-degree-monic-null-poly} and
\ref{theorem:null-poly-omega1-pd-all}). In addition, it is found
that one can enumerate\footnote{In this paper, the word ``enumerate"
means to list all null polynomials modulo $m$, not only to get the
number of all null polynomials modulo $m$. It is obvious that the
former is much more strong than the latter.} all null polynomials
modulo $m$ (see theorems in Sec.
\ref{subsection:null-poly-pd-enumerate-all}, especially Theorem
\ref{theorem:null-poly-enumerate-all}).

The most natural application of a monic null polynomial $f(x)$
modulo $m$ is on the reduction of high-degree congruences modulo
$m$. The null polynomial $f(x)=x^p-x$ has been widely used to reduce
congruences of degrees $>p$ modulo $p$ before solving them
\cite{HW:NumberTheory1979, NZM:NumberTheory1991,
Eynden:NumberTheory1987, Rosen:NumberTheory1993,
Rose:NumberTheory1994, Pan:ConciseNumberTheory1998}. Another
application of the study on null polynomials modulo $m$ is to
distinguish and enumerate equivalent polynomials modulo $m$, i.e.,
the integer polynomials that induce the same ``polynomial functions
modulo $m$" \cite{Mullen:PolyFun-mod:AMH1984}\footnote{In the
algebra literature, there exist ``polynomial functions of the finite
filed $\mathbb{F}_q$" \cite{LauschNobauer:AlgebraPoly1973}. Note
that when $q$ is not a prime, ``polynomial functions modulo $q$" and
``polynomial functions of $\mathbb{F}_q$" are conceptually
different, since $(\mathbb{Z}_q,+,\cdot)$ is just a ring, not a
finite field.}. This is based on a simple theorem (Theorem
\ref{theorem:Equ-Poly-Null-Poly}): \textit{two integer polynomials
$f_1(x)$ and $f_2(x)$ are equivalent modulo $m$ if and only if
$f_1(x)-f_2(x)$ is a null polynomial modulo $m$.} In fact, this
research on null polynomials modulo $m$ was stirred by a pervious
study on enumerating all distinct permutations modulo $m$ induced
from ``permutation polynomials modulo $m$"
\cite{Rivest:PPmod2w:FFTA2001}\footnote{Similarly, in the algebra
literature, there exist ``permutation polynomials of the finite
filed $\mathbb{F}_q$" \cite{LidlMullen:PermPolyI:AMM1988,
LidlMullen:PermPolyII:AMM1990, Mullen:PermPolySurvey:FFaTA1995,
Lidl:FF1997}. As mentioned in previous footnote, when $q$ is a prime
power, ``permutation polynomial modulo $q$" and ``permutation
polynomials of $\mathbb{F}_q$" are conceptually different.}. The
results obtained in this paper can be used to get an exact
estimation of the number of distinct permutations modulo $p^d$
induced from polynomials of degree $n\geq 2p-1$ modulo $p^d$. More
details on the applications of null polynomials for enumerating
permutation polynomials modulo $p^d$ will be given later in a
revised version of [arXiv:math.NT/0509523, 2005]. In addition, we
believe that the study on null polynomials modulo $p^d$ is useful to
reveal some subtle features of the complete systems of polynomial
residues modulo a prime and its powers.

It is well-known that permutation polynomials modulo $m$ can be used
in cryptography and coding \cite{PPinRSA-Crypt83,
Cryptanalysis-PP-PKC-EuroCrypt84, Mullen:PP-NFSR:IEEETIT1989,
PP2Coding2005, RMT:DicksonPoly1993}. Apparently, null polynomials
modulo $m$ can serve as a tool to analyze the security and
performance of the designed ciphers or coding schemes based on
permutation polynomials modulo $m$. For example, the least degree of
all null polynomials modulo $m$ gives an upper bound of the number
of all coefficients of the permutation polynomials used in
cryptography and coding. As shown in Theorem
\ref{theorem:null-poly-omega1-pd-all} of this paper, when $m=p^d$,
the least degree is generally much less than $p^d$, which means that
one has to be more careful when using permutation polynomials to
design cryptosystems.

This paper is organized as follow. Section
\ref{section:Premilinaries} gives some preliminary definitions and
lemmas, as preparations for future discussions. In Sec.
\ref{section:CompositeCase}, we point out that null polynomials
modulo a composite $m$ can be studied via null polynomials modulo
each prime power of $m$. The main body of this paper is Sec.
\ref{section:PrimePowerCase}, in which we discuss null polynomials
modulo a prime and prime powers. In Sec.
\ref{subsection:null-poly-p}, null polynomials modulo $p$ is studied
and it is pointed out that $f(x)=x^p-x$ is the only monic null
polynomial of degree $p$ modulo $p$. In Sec.
\ref{subsection:null-poly-pd-trivial}, some trivial results on null
polynomials modulo $p^d$ ($d\geq 1$) are given. Then, when $2\leq
d\leq p(p+1)+1$, null polynomials modulo $p^d$ are carefully studied
in Sec. \ref{subsection:null-poly-2-I3}, which forms a basis of the
general results for $d\geq 1$ given in Sec.
\ref{subsection:null-poly-pd}. In Sec.
\ref{subsection:null-poly-pd-enumerate-all}, it is studied how to
enumerate all null polynomials modulo $p^d$, based on the results
given in Sec. \ref{subsection:null-poly-pd}.

\section{Preliminaries}
\label{section:Premilinaries}

This section lists a number of definitions and notations used
throughout in this paper. Some preliminary lemmas are also given to
simplify the discussions in this paper. I try to keep the
definitions, notations and lemmas as simple as possible. Please feel
free to contact me if you have some ideas of making them simpler,
more elegant, more beautiful, and/or more rigorous in mathematics.

\subsection{Some Simple Lemmas on Congruences}

The three lemmas will be used in this paper without explicit
citations.

\begin{lemma}[Theorem 2.2 in \cite{NZM:NumberTheory1991}]\label{lemma:polynomial-mod}
Assume $f(x)$ is an integer polynomial. If $x_1\equiv x_2\pmod m$,
then $f(x_1)\equiv f(x_2)\pmod m$.
\end{lemma}

\begin{lemma}
If $a\equiv 0\pmod{m_1}$ and $b\equiv 0\pmod{m_2}$, then $ab\equiv
0\pmod{m_1m_2}$.
\end{lemma}
\begin{proof}
From $a\equiv 0\pmod{m_1}$ and $b\equiv 0\pmod{m_2}$, there exist
$k_1,k_2\in\mathbb{Z}$ such that $a=k_1m_1$ and $b=k_2m_2$. So,
$ab=k_1k_2m_1m_2\equiv 0\pmod{m_1m_2}$.
\end{proof}

The most frequently used form of the above lemma is as follows: if
$f_1(x)$ and $f_2(x)$ are null polynomials modulo $p^{d_1}$ and
$p^{d_2}$, respectively, then $f_1(x)f_2(x)$ is a null polynomial
modulo $p^{d_1d_2}$. See Sec. \ref{subsection:NullPolyDefinition}
for the formal definition of ``null polynomials modulo $m$".

\begin{lemma}
Assume $\bm{A}$ is an $n\times n$ matrix, $\bm{X}$ is a vector of
$n$ unknown integers, and $\bm{B}$ is a vector of $n$ integers. If
$|\bm{A}|$ is relatively prime to $m$, i.e., $\gcd(|\bm{A}|,m)=1$,
then $\bm{A}\bm{X}\equiv\bm{B}\pmod m$ has a unique set of
incongruent solutions
$\bm{X}\equiv\overline{\Delta}(\mathrm{adj}(\bm{A}))\bm{B}\pmod m$,
where $\overline{\Delta}$ is an inverse of $\Delta=|\bm{A}|$ modulo
$m$ and $\mathrm{adj}(\bm{A})$ is the adjoint of $\bm{A}$.
\end{lemma}
\begin{proof}
This lemma is a direct result of Theorem 3.18 in
\cite{Rosen:NumberTheory1993} (see pages 151 and 152).
\end{proof}

\subsection[Polynomial Congruences Modulo $m$]%
{Polynomial Congruences Modulo $\bm{m}$}

The following definition is from Chap. VII of
\cite{HW:NumberTheory1979} and related concepts are slightly
extended.

\begin{definition}
Given two integer polynomials of degree $n$:
$f(x)=a_nx^n+\cdots+a_1x+a_0$ and $g(x)=b_nx^n+\cdots+b_1x+b_0$, if
$\forall i=0\sim n$, $a_i\equiv b_i\pmod m$, we say \uline{$f(x)$ is
congruent to $g(x)$ modulo $m$}, or $f(x)$ and $g(x)$ are
\uline{congruent (polynomials) modulo $m$}, which is denoted by
$f(x)\equiv g(x)\pmod m$. On the other hand, if $\exists
i\in\{1,\cdots,n\}$, such that $a_i\not\equiv b_i\pmod m$, we say
$f(x)$ and $g(x)$ are \uline{incongruent (polynomials) modulo $m$},
denoted by $f(x)\not\equiv g(x)\pmod m$.
\end{definition}

\begin{definition}
A \uline{polynomial congruence (residue) class modulo $m$} is a set
of all polynomials congruent to each other modulo $m$.
\end{definition}

\begin{definition}
A set of polynomials of degree $n$ modulo $m$ is a \uline{complete
system of polynomial residues of degree $n$ modulo $m$}, if for
every polynomial of degree $n$ modulo $m$ there is one and only one
congruent polynomial in this set.
\end{definition}

\begin{lemma}
The following set of polynomials is a complete system of polynomial
residues of degree $n$ modulo $m$:
\[
\mathbb{F}[x]=\left\{f(x)=a_nx^n+\cdots+a_1x+a_0\left|a_n\in\{1,\cdots,m-1\},a_{n-1},\cdots,a_0\in\{0,\cdots,m-1\}\right.\right\}.
\]
\end{lemma}
\begin{proof}
Assume $f(x)=a_nx^n+\cdots+a_1x+a_0$ is a polynomial of degree $n$
modulo $m$. Choose $a_i^*=(a_i\bmod m)\in\{0,\cdots,m-1\}$ ($i=0\sim
n$), then $f^*(x)=a_n^*x^n+\cdots+a_1^*x+a_0^*\in\mathbb{F}$ is
congruent to $f(x)$. Assume that another polynomial
$g(x)=b_nx^n+\cdots+b_1x+b_0\in\mathbb{F}$ is also congruent to
$f(x)$. Then, $\forall i=0\sim n$, $b_i\equiv a_i^*\pmod m$. Since
$\{0,\cdots,m-1\}$ is a complete set of residues modulo $m$,
$b_i=a_i^*$. This means that $g(x)=f^*(x)$. This completes the proof
of this lemma.
\end{proof}

\begin{definition}
A set of polynomials of degree $\leq n$ modulo $m$ is a
\uline{complete system of polynomial residues of degree $\leq n$
modulo $m$}, if for every polynomial of degree $\leq n$ modulo $m$
there is one and only one congruence polynomial.
\end{definition}

\begin{lemma}
The following set of polynomials is a complete system of polynomial
residues of degree $n$ modulo $m$:
\[
\mathbb{F}[x]=\left\{f(x)=a_nx^n+\cdots+a_1x+a_0\left|a_n,a_{n-1},\cdots,a_0\in\{0,\cdots,m-1\}\right.\right\}.
\]
\end{lemma}
\begin{proof}
The proof is similar to the above lemma.
\end{proof}

\subsection[Polynomial Functions Modulo $m$]{Polynomial Functions
Modulo $\bm{m}$}

\begin{definition}
If a function over $\{0,\cdots,m-1\}$ can be represented by a
polynomial modulo $m$, we say this function is \uline{polynomial
modulo $m$}.
\end{definition}

\begin{lemma}
Assume $p$ is a prime. Then, any function over $\{0,\cdots,p-1\}$ is
polynomial modulo $p$.
\end{lemma}
\begin{proof}
Assume $f(x)=a_nx^n+\cdots+a_1x+a_0$ is a polynomial of degree
$n\geq p-1$ modulo $p$. Given a function
$F:\{0,\cdots,p-1\}\to\{0,\cdots,p-1\}$, one has the following
system of congruences:
\[
\left[\begin{matrix}%
1 & 0 & 0^2 & \cdots & 0^{p-1}\\
1 & 1 & 1^2 & \cdots & 1^{p-1}\\
1 & 2 & 2^2 & \cdots & 2^{p-1}\\
\vdots & \vdots & \vdots & \ddots & \vdots\\
1 & p-1 & (p-1)^2 & \cdots & (p-1)^{p-1}
\end{matrix}\right]
\left[\begin{matrix}%
a_0\\
a_1\\
a_2\\
\vdots\\
a_{p-1}
\end{matrix}\right]\equiv
\left[\begin{matrix}%
F(0)\\
F(1)-\sum_{i=p}^na_i\\
F(2)-\sum_{i=p}^n2^ia_i\\
\vdots\\
F(p-1)-\sum_{i=p}^n(p-1)^ia_i
\end{matrix}\right]\pmod{p^d}.
\]
Since the matrix at the left side is a Vondermonde matrix, one can
see its determinant is relatively prime to $p$. So, for each
combination of $a_p,\cdots,a_n$, there is a unique set of
incongruent solutions of $a_0,\cdots,a_{p-1}$. Thus this lemma is
proved.
\end{proof}

\subsection[Equivalent Polynomials Modulo $m$]{Equivalent Polynomials Modulo $\bm{m}$}

The concept of equivalent polynomial modulo $m$ is used to describe
incongruent but equivalent (for any integer) polynomials modulo $m$.

\begin{definition}
Two integer polynomials $f(x)$ and $g(x)$ are \uline{equivalent
(polynomials) modulo $m$} if $\forall x\in\mathbb{Z}$, $f(x)\equiv
g(x)\pmod m$. In other words, two polynomials are equivalent modulo
$m$ if they derive the same polynomial function modulo $m$.
\end{definition}
Note that two equivalent polynomials modulo $m$ may not be congruent
modulo $p$, and may have distinct degrees. As a typical example,
when $p$ is a prime, $f(x)=x^p$ and $g(x)=x$ are equivalent
polynomials modulo $p$.

\begin{lemma}\label{lemma:equivalent-poly-degree1}
Two polynomials of degree 1 modulo $m$, $f(x)=a_1x+a_0$ and
$g(x)=b_1x+b_0$, are equivalent polynomials modulo $m$ if and only
if $f(x)\equiv g(x)\pmod m$, i.e., $a_1\equiv b_1\pmod m$ and
$a_0\equiv b_0\pmod m$.
\end{lemma}
\begin{proof}
The ``if" part is obvious from the definition of equivalent
polynomials modulo $m$, so we focus on the ``only if" part. Since
$f(x)$ and $g(x)$ are equivalent polynomials modulo $m$, then
$\forall x\in\{0,\cdots,m-1\}$,
$f(x)-g(x)=(a_1-b_1)x+(a_0-b_0)\equiv 0\pmod m$. Choosing $x\equiv
0\pmod m$, one has $a_0\equiv b_0\pmod m$. Then, choosing $x\equiv
1\pmod m$, one has $a_1\equiv b_1\pmod m$. Thus this lemma is
proved.
\end{proof}

\begin{lemma}\label{lemma:equivalent-poly-a0}
Two polynomials, $f(x)=a_{n_1}x^{n_1}+\cdots+a_0$ and
$g(x)=b_{n_2}x^{n_2}+\cdots+b_0$, are equivalent polynomials modulo
$m$, then $a_0\equiv b_0\pmod m$.
\end{lemma}
\begin{proof}
Choosing $x=0$, one has $f(x)-g(x)=a_0-b_0\equiv 0\pmod m$. This
lemma is proved.
\end{proof}
\begin{corollary}
Two polynomials, $f(x)=a_nx^n+\cdots+a_2x^2+a_0$ and
$g(x)=a_nx^n+\cdots+a_2x^2+b_0$, are equivalent polynomials modulo
$m$ if and only if $a_0\equiv b_0\pmod m$.
\end{corollary}

\begin{lemma}\label{lemma:equivalent-poly-pd}
Assume $p$ is a prime and $d\geq 1$. Two polynomials,
$f(x)=a_{p-1}x^{p-1}+\cdots+a_0$ and
$g(x)=b_{p-1}x^{p-1}+\cdots+b_0$, are equivalent polynomials modulo
$p^d$ if and only if $f(x)\equiv g(x)\pmod{p^d}$.
\end{lemma}
\begin{proof}
The ``if" part is obvious true, from the definition of equivalent
polynomials modulo $p^d$. So, we focus on the ``only if" part. From
$f(x)-g(x)\equiv 0\pmod{p^d}$, choosing $x=0\sim p-1$, one can get
the following system of congruences in the matrix form
$\bm{A}\bm{X}_{a-b}\equiv \bm{B}\pmod{p^d}$:
\begin{equation}
\left[\begin{matrix}%
1 & 0 & 0^2 & \cdots & 0^{p-1}\\
1 & 1 & 1^2 & \cdots & 1^{p-1}\\
1 & 2 & 2^2 & \cdots & 2^{p-1}\\
\vdots & \vdots & \vdots & \ddots & \vdots\\
1 & p-1 & (p-1)^2 & \cdots & (p-1)^{p-1}
\end{matrix}\right]
\left[\begin{matrix}%
a_0-b_0\\
a_1-b_1\\
a_2-b_2\\
\vdots\\
a_{p-1}-b_{p-1}
\end{matrix}\right]\equiv
\left[\begin{matrix}%
f(0)-g(0)\\
f(1)-g(1)\\
f(2)-g(2)\\
\vdots\\
f(p-1)-g(p-1)
\end{matrix}\right]\equiv
\left[\begin{matrix}%
0\\
0\\
0\\
\vdots\\
0
\end{matrix}\right]\pmod{p^d}.\label{equation:equivalent-poly-pd}
\end{equation}
Since $\bm{A}$ is a Vandermonde sub-matrix, one can get
$|\bm{A}|=\prod_{0\leq i<j\leq p-1}(j-i)$
\cite[\S4.4]{Zhang:MatrixTheory1999}. From $p$ is a prime and
$1\leq(j-i)\leq p-1$, one has $\gcd(|\bm{A}|,p^d)=1$. Thus, the
above system of congruences has a unique set of incongruent
solutions. So, $\forall i=0\sim p-1$, one has $a_i\equiv
b_i\pmod{p^d}$. This completes the proof of this lemma.
\end{proof}
Note that in the above lemma $f(x)$ and $g(x)$ may be polynomials of
degree less than $p-1$ modulo $p^d$. In this case, the matrix at the
left side of the system of congruences may have a smaller size, but
its determinant is still relatively prime to $p^d$.

\begin{corollary}
Assume $p$ is a prime. Two polynomials, $f(x)=a_nx^n+\cdots+a_0$ and
$g(x)=b_nx^n+\cdots+b_0$, are equivalent polynomials modulo $p$ if
and only if $(f(x)\bmod(x^p-x))\equiv(g(x)\bmod(x^p-x))\pmod p$.
\end{corollary}
\begin{proof}
This corollary is a direct result of the above lemma and Fermat's
Little Theorem.
\end{proof}

\subsection[Null Polynomials modulo $m$]{Null Polynomials modulo $\bm{m}$}
\label{subsection:NullPolyDefinition}

\begin{definition}
A polynomial $f(x)$ of degree $n\geq 0$ modulo $m$ is a \uline{null
polynomial of degree $n$ modulo $m$}, if $\forall x\in\mathbb{Z}$,
$f(x)\equiv 0\pmod m$. Specially, $f(x)=0$ is a trivial null
polynomial of degree 0 modulo $m$.
\end{definition}

In the following, we give some simple lemmas on null polynomials
modulo $m$. The proofs of the lemmas are very simple, so they are
omitted here.

\begin{lemma}
If $f(x)=a_nx^n+\cdots+a_1x+a_0$ is a null polynomial modulo $m$,
then $a_0\equiv 0\pmod m$.
\end{lemma}

\begin{lemma}
Given any null polynomial $f(x)$ modulo $m$, $af(x)$ will still be a
null polynomial modulo $m$, where $a$ is an arbitrary integer.
\end{lemma}
\begin{lemma}
A polynomial $f(x)$ is a null polynomial modulo $m$, if and only
$af(x)$ is a null polynomial modulo $m$, where $\gcd(a,m)=1$.
\end{lemma}
\begin{lemma}\label{lemma:null-poly-transitivity}
If $f(x)$ is a null polynomial modulo $m$ and $a\mid m$, then $f(x)$
is still a null polynomial modulo $a$.
\end{lemma}
The most frequently used form of the above lemma is as follows: if
$f(x)$ is a null polynomial modulo $p^d$, then $f(x)$ is still a
null polynomial modulo $p^i$ for any integer $i\leq d$.

\begin{definition}
Denote the least integer $n\geq 1$ such that there exists a null
polynomial of degree $n$ modulo $m$ by $\omega_0(m)$ and call it
\uline{the least null-polynomial degree modulo $m$}. Denote the
least integer $n\geq 1$ such that there exists a \textbf{monic} null
polynomial of degree $n$ modulo $m$ by $\omega_1(m)$ and call it
\uline{the least monic null-polynomial degree modulo $m$}. A (monic)
null polynomial of degree $\omega_0(m)$ or $\omega_1(m)$ is called
\uline{a least-degree (monic) null polynomial modulo $m$}. \iffalse
Sometime we use $\omega(m)$ to denote $\omega_1(m)$, since
$\omega_1(m)$ is much more important than $\omega_0(m)$.\fi
\end{definition}

\begin{lemma}
Every polynomial of degree $\geq\omega_1(m)$ modulo $m$ has one
equivalent polynomial of degree $\leq\omega_1(m)-1$ modulo $m$.
\end{lemma}
\begin{proof}
It is a direct result of the following two fact that there exists a
monic null polynomial of degree $\omega_1(m)$ modulo $m$.
\end{proof}

\begin{definition}
Assume $p$ is a prime and $d\geq 1$. If $f(x)$ is a null polynomial
modulo $p^i$ for any integer $i\leq d$ but not a polynomial modulo
$p^{d+1}$, we say $f(x)$ is a \uline{null polynomial modulo $p$ up
to order $d$} and $d$ is \uline{the order of the null polynomial
$f(x)$ modulo $p$}. Note that $f(x)$ is a null polynomial modulo
$p^d$ can ensure that $f(x)$ is a null polynomial modulo $p^i$ for
any integer $i\leq d$.
\end{definition}

\section{Null Polynomial modulo $\bm{m=\prod_{i=1}^rp_i^{d_i}}$}
\label{section:CompositeCase}

\begin{theorem}\label{theorem:null-poly-composite}
Assume $p_1$, $\cdots$, $p_r$ are $r$ distinct prime numbers and
$d_1$, $\cdots$, $d_r\geq 1$. A polynomial $f(x)$ is a null
polynomial modulo $m=\prod_{i=1}^rp_i^{d_i}$, if and only if
$\forall i=1\sim r$, $f(x)$ is a null polynomial modulo $p_i^{d_i}$.
\end{theorem}
\begin{proof}
It is a direct result of the Chinese remainder theorem.
\end{proof}

\begin{theorem}
Assume $p_1$, $\cdots$, $p_r$ are $r$ distinct prime numbers, $d_1$,
$\cdots$, $d_r\geq 1$ and $m=\prod_{i=1}^rp_i^{d_i}$. If $f_1(x)$,
$\cdots$, $f_r(x)$ are null polynomials of degree $n_1$, $\cdots$,
$n_r$ modulo $p_1^{d_1}$, $\cdots$, $p_r^{d_r}$, respectively, then
there exists one and only one null polynomial $f(x)$ of degree
$D=\max_{i=1}^r(d_i)$ modulo $m$ in each complete system of
polynomial residues modulo $m$, such that $f(x)\equiv
f_i(x)\pmod{p_i^{d_i}}$ holds for $i\in\{1,\cdots,r\}$.
\end{theorem}
\begin{proof}
Applying the Chinese remainder theorem on each coefficient of the
$r$ polynomials\footnote{For polynomials of degree less than $D$,
the higher coefficients are assigned to be zeros.}, one can
immediately prove this theorem.
\end{proof}

\begin{theorem}\label{theorem:null-poly-omega-composite}
Assume $p_1$, $\cdots$, $p_r$ are $r$ distinct prime numbers, $d_1$,
$\cdots$, $d_r\geq 1$ and $m=\prod_{i=1}^rp_i^{d_i}$. Then,
$\omega(m)=\max_{i=1}^r(\omega(p_i^{d_i}))$, i.e., the least degree
of null polynomials modulo $m$ is $\max_{i=1}^r(\omega(p_i^{d_i}))$.
Here, $\omega(p_i^{d_i})$ can be either $\omega_0(p_i^{d_i})$ or
$\omega_1(p_i^{d_i})$.
\end{theorem}
\begin{proof}
From the above theorem, one can find a (monic) null polynomial of
degree $D$ modulo $m$, so $\omega(m)\leq D$. Next, assume there
exists another (monic) null polynomial $g(x)$ of degree $\leq D-1$
modulo $m$. Then, $g(x)$ is also a (monic) null polynomial modulo
each $p_i^{d_i}$. This means $\max_{i=1}^r(\omega(p_i^{d_i}))\leq
D-1$. We get a contradiction, so $\omega(m)=D$. This theorem is thus
proved.
\end{proof}

With the above theorem, the composite case can be handled easily by
handling the $r$ prime power cases.

\begin{corollary}
Assume $p_1$, $\cdots$, $p_r$ are $r$ distinct prime numbers, $d_1$,
$\cdots$, $d_r\geq 1$, $m=\prod_{i=1}^rp_i^{d_i}$ and
$D=\max_{i=1}^r(\omega_1(p_i^{d_i}))$. For each polynomial $f(x)$ of
degree $\geq D$ modulo $m$, there exists an equivalent polynomial of
degree $\leq D-1$ modulo $m$.
\end{corollary}
\begin{proof}
This corollary is a direct result of the above theorem.
\end{proof}

\section{Null Polynomials modulo $\bm{p^d}$ ($\bm{d\geq 1}$)}
\label{section:PrimePowerCase}

In this section, the following questions are focused.

\begin{enumerate}
\item
What are the values of $\omega_0(p^d)$ and $\omega_1(p^d)$, i.e.,
the smallest integer $n$ such that there exists at least one (monic)
null polynomial modulo $p^d$?

\item
What is the number of (monic) null polynomials of degree $n$ modulo
$p^d$?

\item
Is it possible to enumerate all incongruent (monic) null polynomials
of degree $n$ modulo $p^d$?
\end{enumerate}

\subsection[Null Polynomials modulo $p$]{Null Polynomials modulo $\bm{p}$}
\label{subsection:null-poly-p}

\begin{theorem}\label{theorem:null-poly-omega-p}
Assume $p$ is a prime. Then, $\omega_0(p)=\omega_1(p)=p$.
\end{theorem}
\begin{proof}
It is a direct result of Lemma \ref{lemma:equivalent-poly-pd}, since
any two equivalent polynomials of degree $\leq p-1$ modulo $p$ are
congruent modulo $p$ and $f(x)=x^p-x$ is a null polynomial modulo
$p$.
\end{proof}

\begin{theorem}\label{theorem:null-poly-p}
Assume $p$ is a prime and $F(x)$ is a null polynomial of degree $p$
modulo $p$. Then, $f(x)$ is a null polynomial modulo $p$ if and only
if $f(x)\equiv F(x)q(x)\pmod p$, where $q(x)$ is an arbitrary
polynomial modulo $p$.
\end{theorem}
\begin{proof}
The ``if" part is obviously true. Let us prove the ``only if" part.
Dividing $f(x)$ by $F(x)$, one can get $f(x)=F(x)q(x)+r(x)$, where
$r(x)$ is a polynomial of degree $\leq p-1$ modulo $p$. Since
$F(x)q(x)$ is a null polynomial modulo $p$, $r(x)$ is also a null
polynomial modulo $p$. From Lemma \ref{lemma:equivalent-poly-pd},
$r(x)$ is congruent to zero polynomial modulo $p$. Thus, $f(x)\equiv
F(x)q(x)\pmod p$. Thus, this theorem is proved.
\end{proof}

From the above theorem, one can enumerate all null polynomials of
degree $n\geq p$ modulo $p$.

\begin{corollary}
Assume $p$ is a prime. Then, $f(x)$ is a null polynomial modulo $p$
if and only if $f(x)\equiv (x^p-x)q(x)\pmod p$, where $q(x)$ is an
arbitrary polynomial modulo $p$.
\end{corollary}
\begin{proof}
It is a direct result of the above theorem, since $x^p-x$ is a null
polynomial modulo $p$.
\end{proof}
\begin{corollary}
Assume $p$ is a prime and $f(x)$ is a null polynomial of degree $p$
modulo $p$, then $f(x)\equiv a(x^p-x)\pmod p$, where $\gcd(a,p)=1$.
That is, $x^p-x$ is the \textbf{only one monic} null polynomial
modulo $p$.
\end{corollary}
\begin{proof}
It is a direct result of the above theorem. Note that $f(x)$
congruent to zero polynomial modulo $p$ if $\gcd(a,p)>1$, i.e.,
$a\equiv 0\pmod p$.
\end{proof}
\begin{corollary}
Assume $p$ is a prime, then $x(x-1)\cdots(x-(p-1))\equiv(x^p-x)\pmod
p$.
\end{corollary}

\begin{remark}
Note that the above corollary is generally proved in number theory
via Lagrange's Theorem. Since Lemma \ref{lemma:equivalent-poly-pd}
and Theorem \ref{theorem:null-poly-omega-p} do not depend on the
Lagrange's Theorem, so we give an independent proof of the
well-known result. Note that Wilson's Theorem can be derived from
this corollary by choosing $x=0$.
\end{remark}

\newcommand\NPp[1][p]{\mathcal{F}_{#1}(x)}
\newcommand\NPpf[1][p]{\mathcal{F}_{#1}}
\newcommand\NPpd[2][p]{\mathcal{F}_{#1,#2}}
\newcommand\degree{\mathrm{deg}}

\begin{definition}
In the following of this paper, to facilitate the discussion, define
$\NPp=\prod_{i=0}^{p-1}(x-i)=x(x-1)\cdots(x-(p-1))$. This special
polynomial will be frequently used to derive some important results.
\end{definition}

\begin{definition}
Assume $p$ is a prime. Define $\Lambda(x)=\prod_{0\leq j\leq p-1
\atop j\not\equiv x\pmod p}(x-j)$.
\end{definition}

\begin{lemma}
Assume $p$ is a prime, then $\forall x\in\mathbb{Z}$,
$\Lambda(x)\equiv(p-1)!\equiv -1\pmod p$.
\end{lemma}
\begin{proof}
Assuming $i=(x\bmod p)\in\{0,\cdots,p-1\}$, one has
$\Lambda(x)\equiv\Lambda(i)=\prod_{0\leq j\leq p-1 \atop j\neq
i}(i-j)\equiv\prod_{i+1\leq j\leq p-1}(p+i-j)\prod_{0\leq j\leq
i-1}(i-j)=(p-1)\cdots(i+1)i!=(p-1)!\equiv -1\pmod p$. In fact, this
lemma is true since $\{x-j\}_{0\leq j\leq p-1 \atop j\not\equiv
x\pmod p}$ actually forms a reduced system of residues modulo $p$.
\end{proof}

\begin{theorem}\label{theorem:base-null-poly-p}
Assume $p$ is a prime. Then, $\forall i\in\mathbb{Z}$ and $\forall
j\in\{0,\cdots,p-1\}$, $\NPpf(ip+j)\equiv-ip\pmod{p^2}$, i.e.,
$\frac{\NPpf(ip+j)}{p}\equiv-i\pmod p$.
\end{theorem}
\begin{proof}
From the above lemma, one has $\frac{\NPpf(ip+j)}{p}=
\frac{(ip)\Lambda(ip+j)}{p}=i\Lambda(ip+j)\equiv-i\pmod p$.
\end{proof}
Yet another form of the above theorem is as follows.
\begin{theorem}\label{theorem:base-null-poly-p-form2}
Assume $p$ is a prime. Then, $\forall x\in\mathbb{Z}$,
$\frac{\NPpf(x)}{p}=\lfloor x/p\rfloor\Lambda(x)\equiv-\lfloor
x/p\rfloor\pmod p$.
\end{theorem}
The above theorem is very important as a basic feature of $\NPp$ to
generalize the results modulo $p$ to modulo any power of $p$.

\subsection[Null Polynomials modulo $p^d$ ($d\geq 1$): Some Trivial Results]%
{Null Polynomials modulo $\bm{p^d}$ ($\bm{d\geq 1}$): Some Trivial
Results} \label{subsection:null-poly-pd-trivial}

This subsection gives some trivial results on null polynomials
modulo $p^d$, some of which (especially Corollary
\ref{corollary:null-poly-hierarchy-pd}) will be frequently cited
later to get some more important results.

\begin{theorem}
Assume $p$ is a prime and $d\geq 1$. Then, $\omega_0(p^d)=p$.
\end{theorem}
\begin{proof}
Note that $p^{d-1}(x^p-x)$ is a null polynomial modulo $p^d$, so
$\omega_0(p^d)=p$.
\end{proof}

From the above theorem, one can see that $\omega_0(p^d)$ is trivial
for studying null polynomial modulo $p^d$. So in following we will
focus on $\omega_1(p^d)$ only. At first, we introduce some
preliminary lemmas for further discussions. They will be cited later
without explicit citations.

\begin{lemma}
Assume $p$ and $d\geq 1$. If $f(x)=a_nx^n+\cdots+a_1x$ is a null
polynomial of degree $<\omega_1(p^d)$ modulo $p^d$, then one of the
following results holds:
\begin{itemize}
\item
$f(x)\equiv 0\pmod{p^d}$, i.e., $f(x)=p^df^*(x)$, where $f^*(x)$ is
any integer polynomial;

\item
$f(x)\equiv ap^if^*(x)\pmod{p^d}$, where $\gcd(a,p)=1$,
$i\in\{1,\cdots,d\}$ and $f^*(x)$ is a monic null polynomial modulo
$p^{d-i}$.
\end{itemize}
\end{lemma}
\begin{proof}
When $f(x)\equiv 0\pmod{p^d}$, one has $f(x)=p^df^*(x)$, where
$f^*(x)$ can be any polynomials. When $f(x)\not\equiv 0\pmod{p^d}$,
assuming $p\nmid a_n$, then $\gcd(a_n,p)=1$, so there exists an
inverse $\bar{a}_n$ such that $\bar{a}_na_n\equiv 1\pmod{p^d}$.
Multiplying $f(x)$ by $\bar{a}_n$, one gets a monic null polynomial
$f^*(x)=x^n+\cdots+\bar{a}_na_1x$ modulo $p^d$. This conflicts with
the fact $n<\omega_1(p^d)$. So $p\mid a$ is always true and one has
$a_n=ap^i$, where $\gcd(a,p)=1$ and $i\in\{1,\cdots,d-1\}$. Choosing
$f^*(x)=$ This proves this lemma.
\end{proof}
\begin{corollary}
Assume $p$ and $d\geq 1$. If $f(x)=a_nx^n+\cdots+a_1x$ is a null
polynomial of degree $<\omega_1(p^d)$ modulo $p^d$, then
$f(x)=pf^*(x)\pmod{p^d}$, where $f^*(x)$ is a null polynomial modulo
$p^{d-1}$.
\end{corollary}
\begin{proof}
It is a direct result of the above lemma.
\end{proof}

\begin{lemma}\label{lemma:null-poly-hierarchy-pd}
Assume $p$ is a polynomial, $1\leq d_1<\cdots<d_n<d$ and
$F_1(x),\cdots,F_n(x)$ are monic null polynomials of degree
$\omega_1(p^{d_1})<\cdots<\omega_1(p^{d_n})$ modulo
$p^{d_1},\cdots,p^{d_n}$, respectively. If $f(x)$ is a null
polynomial modulo $p^d$, then
\[
f(x)\equiv\sum_{i=n}^1F_i(x)p^{n-i}q_i(x)+p^nq_0(x)\pmod{p^d},
\]
where $q_0(x)$ is a polynomial of degree less than
$\omega_1(p^{d_1})$ modulo $p^d$ and $q_i(x)$ ($1\leq i\leq n-1$) is
a polynomial of degree less than
$\omega_1(p^{d_{i+1}})-\omega_1(p^{d_i})$ modulo $p^d$.
\end{lemma}
\begin{proof}
We use induction on $n$ to prove this lemma.

When $n=1$, dividing $f(x)$ by $F_1(x)$, one has
$f(x)=F_1(x)q_1(x)+r_1(x)$, where $r_1(x)$ is of degree less than
$\omega_1(p^{d_1})$ modulo $p^d$. Since $f(x)$ and $F_1(x)q_1(x)$
are both null polynomials modulo $p^{d_1}$, $r_1(x)$ is also a null
polynomial modulo $p^{d_1}$. Considering the degree of $r_1(x)$ is
less than $\omega_1(p^{d_1})$, so one has $r_1(x)\equiv
pr_1^*(x)\pmod{p^{d_1}}$ and then $r_1(x)\equiv
pr_1^{(1)}(x)+p^{d_1}r_1^{(1)}(x)\equiv pq_0(x)\pmod{p^d}$. This
proves the case of $n=1$.

Assume this lemma is true for any integer $\leq n-1$, let us prove
the case of $n\geq 2$. Dividing $f(x)$ by $F_n(x)$, one has
$f(x)=F_n(x)q_n(x)+r_n(x)$, where $r_n(x)$ is a polynomial of degree
$<\omega_1(p^{d_n})$ modulo $p^d$. Since $f(x)$ and $F_n(x)q_n(x)$
are both null polynomial modulo $p^{d_n}$, $r_n(x)$ is also a null
polynomial modulo $p^{d_n}$. Considering $r_n(x)$ is of degree less
than $\omega_1(p^{d_n})$, one has $r_n(x)\equiv
pr_n^*(x)\pmod{p^{d_n}}$ and then $r_n(x)\equiv
pr_n^{(1)}(x)+p^{d_n}r_n^{(2)}(x)\equiv pr_n^{**}(x)\pmod{p^d}$,
where $r_n^{**}(x)$ is a null polynomial of degree less than
$\omega_1(p^{d_n})$ modulo $p^{d_{n-1}}$. Then, applying the
hypothesis on $r_n^{**}(x)$, the case of $d\geq 2$ is proved.
\end{proof}

\begin{corollary}\label{corollary:null-poly-hierarchy-pd}
Assume $p$ is a polynomial, $d\geq 2$ and $F_1(x)$, $\cdots$,
$F_{d-1}(x)$ are monic null polynomials of degree $\omega_1(p)$,
$\cdots$, $\omega_1(p^{d-1})$ modulo $p$, $\cdots$, $p^{d-1}$,
respectively. If $f(x)$ is a null polynomial modulo $p^d$, then
\[
f(x)\equiv\sum_{n=d-1}^1F_n(x)p^{d-1-n}q_n(x)+p^{d-1}q_0(x)\pmod{p^d},
\]
where $q_0(x)$ is a polynomial of degree less than $\omega_1(p)$
modulo $p^d$ and $q_n(x)$ ($1\leq n\leq d-2$) is a polynomial of
degree less than $\omega_1(p^{n+1})-\omega_1(p^n)$ modulo $p^d$. If
there exists $1\leq n_1<n_2\leq d-1$ such that
$\omega_1(p^{n_1})=\cdots=\omega_1(p^{n_2})$, then
$q_{n_1}(x)=\cdots=q_{n_2-1}(x)=0$.
\end{corollary}

\begin{explain}
The above corollary is a direct result of Lemma
\ref{lemma:null-poly-hierarchy-pd}. This corollary makes it possible
to use induction on $d$ to derive null polynomials modulo $p^d$ from
null polynomials modulo lower prime powers $p,\cdots,p^{d-1}$.
\end{explain}

\begin{lemma}
Assume $p$ is a prime and $d\geq 1$. Then,
$0\leq\omega_1(p^{d+1})-\omega_1(p^d)\leq p$.
\end{lemma}
\begin{proof}
From Lemma \ref{lemma:null-poly-transitivity}, one can easily get
$\omega_1(p^{d+1})\geq\omega_1(p^d)$. Then, one has
$\omega_1(p^{d+1})\leq\omega_1(p^d)+p$ from the following fact: if
$f(x)$ is a monic null polynomial modulo $p^d$, then $f(x)\NPp$ is a
monic null polynomial modulo $p^{d+1}$. So, one has
$0\leq\omega_1(p^{d+1})-\omega_1(p^d)\leq p$.
\end{proof}

\subsection[Null Polynomials modulo $p^d$: The Case of $2\leq d\leq p(p+1)+1$]%
{Null Polynomials modulo $\bm{p^d}$: The Case of $\bm{2\leq d\leq
p(p+1)+1}$} \label{subsection:null-poly-2-I3}

In this subsection, we study the case of $2\leq d\leq p(p+1)+1$. The
results obtained on these special cases lead to a recursive way to
handle the general case of $d\geq 1$ (as shown in next subsection).

\subsubsection[The Case of $2\leq d\leq p$]{The Case of $\bm{2\leq d\leq p}$}

\begin{lemma}
Assume $p$ is a prime. A polynomial $f(x)$ is a null polynomial
modulo $p^2$, if and only if $f(x)\equiv\NPp
q_1(x)\equiv(\NPp)^2q_1^*(x)+p\NPp q_0^*(x)\pmod{p^2}$, where
$q_1(x)$ is a null polynomial modulo $p$ and $q_0^*(x),q_1^*(x)$ can
be any integer polynomials.
\end{lemma}
\begin{proof}
The ``if" part is obviously true, so we focus on the ``only if"
part.

From Corollary \ref{corollary:null-poly-hierarchy-pd}, one has
$f(x)\equiv \NPp q_1(x)+pq_0(x)\pmod{p^2}$, where $q_0(x)$ is of
degree $\leq p-1$ modulo $p^2$. Choosing $x=0$, one has
$pq_0(0)\equiv 0\pmod {p^2}\Rightarrow q_0(0)\equiv 0\pmod p$.
Choosing $x=ip$, where $i\in\{1,\cdots,p-1\}$, one has
$ip\Lambda(ip)q_1(ip)+pq_0(ip)\equiv 0\pmod {p^2}\Rightarrow
i\Lambda(ip)q_1(ip)+q_0(ip)\equiv 0\pmod p\Rightarrow -iq_1(0)\equiv
0\pmod p\Rightarrow q_1(0)\equiv 0\pmod p$. Next, choosing
$x=j\in\{1,\cdots,p-1\}$, one has $pq_0(j)\equiv
0\pmod{p^2}\Rightarrow q_0(j)\equiv 0\pmod p$. Then, choosing
$x=p+j\in\{p+1,\cdots,2p-1\}$, one has
$p\Lambda(p+j)q_1(p+j)+pq_0(p+j)\equiv 0\pmod{p^2}\Rightarrow
\Lambda(p+j)q_1(j)\equiv 0\pmod p\Rightarrow q_1(j)\equiv 0\pmod p$.
Combining the above results, one can see that both $q_1(x)$ and
$q_0(x)$ are null polynomials modulo $p$. However, since $q_0(x)$ is
a null polynomial of degree less than $p$, one immediately gets
$q_0(x)\equiv 0\pmod p$ and then $pq_0(x)\equiv 0\pmod{p^2}$. Thus,
$f(x)\equiv \NPp q_1(x)\pmod{p^2}$. Considering $q_1(x)\equiv \NPp
q_1^*(x)\pmod p$, one has $q_1(x)\equiv \NPp
q_1^*(x)+pq_0^*(x)\pmod{p^2}$. As a final result, $f(x)\equiv
(\NPp)^2q_1^*(x)+p\NPp q_0^*(x)\pmod{p^2}$, where $q_0^*(x)$ and
$q_1^*(x)$ can be any integer polynomials (without the limit on the
degree modulo $p^2$). This completes the proof of this lemma.
\end{proof}

\begin{remark}
In the above lemma, the two different representations of the
necessary and sufficient conditions of $f(x)$ may have different
usages. The first one, $f(x)\equiv\NPp q_1(x)\pmod{p^2}$, is better
to show the law basic law behind the result and to organize the
proof; while the second one, $f(x)\equiv(\NPp)^2q_1^*(x)+p\NPp
q_0^*(x)\pmod{p^2}$, is better to enumerate all null polynomials
modulo $p^2$. In the following of this section, we continue to adopt
the two representations simultaneously.
\end{remark}

\begin{lemma}
Assume $p$ is a prime, then $\omega_1(p^2)=2p$.
\end{lemma}
\begin{proof}
From the above lemma, to get a monic null polynomial modulo $p^2$,
it is obvious that $q_1(x)\not\equiv 0\pmod p$. Since the least
degrees of $q_1(x)$ and $\NPp$ are both $p$, the least degree of
$f(x)$ is $p+p=2p$, i.e., $\omega_1(p^2)=2p$, where note that
$f_1(x)$ is a monic polynomial modulo $p^2$ if $q_1(x)$ is a monic
polynomial modulo $p^2$. This completes the proof of this lemma.
\end{proof}

\begin{theorem}\label{theorem:null-poly-2-p}
Assume $p$ is a prime and $2\leq d\leq p$. A polynomial $f(x)$ is a
null polynomial modulo $p^d$ if and only if $f(x)\equiv
\sum_{n=d}^1p^{d-n}(\NPp)^nq_n^*(x)\pmod{p^d}$, where $q_1^*(x)$,
$\cdots$, $q_{d-1}^*(x)$ are \textbf{any} polynomials of degree less
than $p$ modulo $p^d$ and $q_d(x)$ is \textbf{any} polynomial of
\textbf{any} degree modulo $p$.
\end{theorem}
\begin{proof}
The ``if" part is obviously true, so we only prove the ``only if"
part via induction on $d$.

When $d=2$, this theorem has been proved above. Let us use induction
on $d$ to prove the case of $3\leq d\leq p$, under the assumption
that this theorem is true for all integers not greater than $d-1$.

Apparently, the hypothesis means that $\omega_1(p^c)=pc$ when $c<d$.
Then, from Corollary \ref{corollary:null-poly-hierarchy-pd}, one has
$f(x)\equiv \sum_{n=d-1}^0p^{d-1-n}(\NPp)^nq_n(x)\pmod{p^d}$, where
$q_0(x)$, $\cdots$, $q_{d-2}(x)$ are of degree less than
$\omega_1(p^{n+1})-\omega_1(p^n)=p$.

Assuming $x=ip+j$, where $i$ is any integer and
$j\in\{0,\cdots,p-1\}$. One can see that $\{ip+j\}$ forms a
completes system of residues modulo $p^d$ when $i$ runs through a
complete system of residues modulo $p^{d-1}$. Substituting $x=ip+j$
into $f(x)\equiv 0\pmod{p^d}$, one has
$\sum_{n=d-1}^0p^{d-1-n}(ip\Lambda(ip+j))^nq_n(ip+j)\equiv
0\pmod{p^d}$. This congruence can be further simplified as
$\sum_{n=d-1}^0(-i)^nq_n(j)\equiv 0\pmod p$, where note that
$\Lambda(ip+j)\equiv -1\pmod p$.

Choosing $i\equiv 0\pmod p$, one gets $q_0(j)\equiv 0\pmod p$ for
any $j$. Considering the degree of $q_0(x)$ is less than $p$,
$q_0(x)\equiv 0\pmod p\Rightarrow p^{d-1}q_0(x)\equiv 0\pmod{p^d}$,
so this term can be removed. The congruence is simplified to be
$\sum_{n=d-1}^1(-i)^nq_n(j)\equiv 0\pmod p$. To solve the value of
each $q_n(j)$ when $n\geq 1$, consider the polynomial
$h_j(x)=\sum_{n=d-1}^1(-1)^nq_n(j)x^n$ of degree $d-1\leq p-1$.
Since $h_j(x)\equiv 0\pmod p$ holds for any integer $x$ and
$j\in\{0,\cdots,p-1\}$, one immediately has $(-1)^nq_n(j)\equiv
0\pmod p$ and then $q_n(j)\equiv 0\pmod p$ for all values of $n\geq
1$ and $j\in\{0,\cdots,p-1\}$. Thus, finally we get the result that
$q_1(x),\cdots,q_{d-1}(x)$ are all null polynomials modulo $p$.
Since the first $d-2$ polynomials are of degree less than $p$, one
immediately has $q_n(x)\equiv pq_n^*(x)\pmod{p^d}$. The last
polynomial is of any degree, so $q_{d-1}(x)\equiv\NPp
q_d^*(x)+pq_{d-1}^*(x)\pmod{p^d}$, where $q_{d-1}^*(x)$ is also of
degree less than $p$ modulo $p^d$. Substituting the $d-1$
congruences into $f(x)$, one has $f(x)\equiv
\sum_{n=d}^1p^{d-n}(\NPp)^nq_n^*(x)\pmod{p^d}$. Thus, the case of
$3\leq d\leq p$ is proved and the proof of this theorem is also
completed.
\end{proof}

\begin{theorem}\label{theorem:null-poly-omega-22p}
Assume $p$ is a prime and $2\leq d\leq p$. Then, $\omega_1(p^d)=pd$.
\end{theorem}
\begin{proof}
It is a direct result of the above theorem, since the least degree
of $q_d^*(x)$ corresponding to a monic null polynomial modulo $p^d$
is 0.
\end{proof}

\subsubsection[The Case of $d=p+1$]{The Case of $\bm{d=p+1}$}

\begin{theorem}\label{theorem:null-poly-p+1-a}
Assume $p$ is a prime. A polynomial $f(x)$ is a null polynomial
modulo $p^{p+1}$, if and only if the following conditions hold
simultaneously:
\begin{itemize}
\item
$f(x)\equiv\sum_{n=p}^1p^{p-n}(\NPp)^nq_n(x)\pmod{p^{p+1}}$, where
$q_1(x)$, $\cdots$, $q_{p-1}(x)$ are polynomial of degree less than
$p$ modulo $p^{p+1}$;

\item
$q_p(x)$ is any polynomial of any degree modulo $p^{p+1}$ and
$\sum_{n=p-1}^1(-i)^nq_n(j)\equiv iq_p(j)\pmod p$ holds for
$i=\lfloor x/p\rfloor$ and $j=x\bmod p$.
\end{itemize}
\end{theorem}
\begin{proof}
At first, let us prove the ``only if" part.

Since $\omega_1(p^p)=p^2$, one has
$f(x)\equiv\sum_{n=p}^0p^{p-n}(\NPp)^nq_n(x)\pmod{p^{p+1}}$.
Assuming $x=ip+j$, where $i$ runs through a complete system of
residues modulo $p^p$ and $j\in\{0,\cdots,p-1\}$. Substituting
$x=ip+j$ into $f(x)\equiv 0\pmod{p^{p+1}}$, one has
$\sum_{n=p}^0p^{p-n}(ip\Lambda(ip+j))^nq_n(ip+j)\equiv
0\pmod{p^{p+1}}$. This congruence can be reduced to be
$\sum_{n=p}^0(-i)^nq_n(j)\equiv 0\pmod p$, where note that
$\Lambda(ip+j)\equiv -1\pmod p$.

Choosing $i\equiv 0\pmod p$, one immediately gets $q_0(j)\equiv
0\pmod p$, which means that $p^pq_0(x)\equiv 0\pmod{p^{p+1}}$, so
this term can be removed. Next, for a given value of $j$, to solve
each of other $q_n(j)$, one has the following system of congruences:
for $i\equiv 1\sim p-1\pmod p$, $\sum_{n=p}^1(-i)^nq_n(j)\equiv
0\pmod p$. Write the $p-1$ congruences as the matrix form after
moving $(-i)^pq_p(j)$ to the right side as follows,
\[
\left[\begin{matrix}%
1 & 1 & \cdots & 1\\
2 & 2^2 & \cdots & 2^{p-1}\\
\vdots & \vdots &\ddots & \vdots\\
(p-1) & (p-1)^2 & \cdots & (p-1)^{p-1}
\end{matrix}\right]
\left[\begin{matrix}%
-q_1(j)\\(-1)^2q_2(j)\\\vdots\\(-1)^{p-1}q_{p-1}(j)
\end{matrix}\right]\equiv
\left[\begin{matrix}%
-(-1)^pq_p(j)\\-(-2)^pq_p(j)\\\vdots\\-(-(p-1))^pq_p(j)
\end{matrix}\right]\equiv
\left[\begin{matrix}%
q_p(j)\\2q_p(j)\\\vdots\\(p-1)q_p(j)
\end{matrix}\right]\pmod p.
\]
Since the matrix at left side is a Vondermande matrix, one can see
that its determinant is relatively prime to $p$. So, for any value
of $q_p(j)$, there exists a unique set of incongruent solutions of
$\{(-1)^nq_n(j)\}_{1\leq n\leq p-1}$ modulo $p$. That is, there is a
unique set of incongruent solutions to $\{q_n(j)\}_{1\leq n\leq
p-1}$ modulo $p$. Since every function over $\{0,\cdots,p-1\}$
corresponds to a unique polynomial of degree $\leq p-1$ modulo $p$,
there exists a unique polynomial $q_n(x)$ ($1\leq n\leq p-1$) of
degree $\leq p-1$ modulo $p$, for each distinct polynomial $q_p(x)$
modulo $p$.

The ``if" part can be easily proved since
$\sum_{n=p-1}^1(-i)^nq_n(j)\equiv iq_p(j)\pmod p$ actually means
$f(x)\equiv 0\pmod{p^{p+1}}$. Thus, this theorem is proved.
\end{proof}

Yet another form of the above theorem is as follow, in which the
free polynomial becomes $q_1(x)$.

\begin{theorem}\label{theorem:null-poly-p+1-b}
Assume $p$ is a prime. A polynomial $f(x)$ is a null polynomial
modulo $p^{p+1}$, if and only if the following conditions hold
simultaneously:
\begin{itemize}
\item
$f(x)\equiv\sum_{n=p}^1p^{p-n}(\NPp)^nq_n(x)\pmod{p^{p+1}}$, where
$q_1(x)$, $\cdots$, $q_{p-1}(x)$ are polynomial of degree less than
$p$ modulo $p^{p+1}$;

\item
$q_1(x)$ is any polynomial of degree less than $p$ modulo $p^{p+1}$
and $\sum_{n=p}^2(-i)^nq_n(j)\equiv -iq_1(j)\pmod p$ holds for
$i=\lfloor x/p\rfloor$ and $j=x\bmod p$.
\end{itemize}
\end{theorem}
\begin{proof}
In the proof of the above theorem, moving $-q_1(j)$ to the right
side of the congruence, one has
\begin{equation}
\left[\begin{matrix}%
1 & 1 & \cdots & 1\\
2^2 & 2^3 & \cdots & 2^p\\
\vdots & \vdots &\ddots & \vdots\\
(p-1)^2 & (p-1)^3 & \cdots & (p-1)^p
\end{matrix}\right]
\left[\begin{matrix}%
(-1)^2q_2(j)\\(-1)^3q_3(j)\\\vdots\\(-1)^pq_p(j)
\end{matrix}\right]\equiv
\left[\begin{matrix}%
-q_1(j)\\-2q_1(j)\\\vdots\\-(p-1)q_1(j)
\end{matrix}\right]\pmod p.\label{equation:null-poly-p+1-b}
\end{equation}
It is obvious that the determinant of the matrix at the left side is
still relatively prime to $p$, so there is a unique solution to each
value of $q_1(j)$. Thus, this theorem is true.
\end{proof}

\begin{theorem}
In the above theorems, 1) $q_1(x)$ or $q_p(x)$ uniquely determines
the polynomial functions derived from all other $p-1$ polynomials
modulo $p$; 2) if $q_1(x)$ or $q_p(x)$ is of degree 0, all other
$p-1$ polynomials are of degree 0; 3) $q_1(j)\equiv 0\pmod
p\Leftrightarrow q_p(j)\equiv 0\pmod p$ and $q_1(j)\not\equiv 0\pmod
p\Leftrightarrow q_p(j)\not\equiv 0\pmod p$.
\end{theorem}
\begin{proof}
The proof of this theorem is actually included in the proof of the
above theorem.
\end{proof}

\begin{remark}
Specially, when $q_1(x)$ or $q_p(x)$ is a null polynomial modulo
$p$, all polynomials are null polynomials modulo $p$. In this case,
$f(x)\equiv\sum_{n=p+1}^1p^{p+1-n}(\NPp)^nq_n^*(x)\pmod{p^{p+1}}$,
where $q_1^*(x),\cdots,q_p^*(x)$ are any polynomials of degree less
than $p$ modulo $p^{p+1}$ and $q_{p+1}^*(x)$ is any polynomial of
any degree modulo $p^{p+1}$. However, null polynomials in this form
are not least-degree null polynomials modulo $p^{p+1}$.
\end{remark}

\begin{theorem}
Assume $p$ is a prime, then $\omega_1(p^{p+1})=p^2$.
\end{theorem}
\begin{proof}
From the above theorem, to get a monic null polynomial modulo
$p^{p+1}$, it is obvious that $q_p(x)\not\equiv 0\pmod p$. Since the
least degree of $q_p(x)$ is $0$ and the degree of $(\NPp)^p$ is
$p^2$, the least degree of $f(x)$ is also $p^2$, i.e.,
$\omega_1(p^{p+1})=p^2$, where note that $f(x)$ is a monic
polynomial modulo $p^3$ if $q_p(x)$ is a monic polynomial modulo
$p^{p+1}$. This theorem is thus proved.
\end{proof}

\begin{example}
When $p=2$, find a monic null polynomial of degree
$\omega_1(2^3)=2^2=4$ modulo $2^3$.
\end{example}
\begin{solution}
In this case, only one congruence is left: $q_1(j)\equiv
-q_2(j)\equiv q_2(j)\pmod p$. Choosing $q_1(x)=q_2(x)=1$, one gets a
monic null polynomial of degree 4 modulo $2^3$:
$f(x)=\uwave{(x^2-x)^2+2(x^2-x)}=\uwave{x^4-2x^3+3x^2-2x}$.
Experiments have been made to verify this null polynomial.
\end{solution}
\begin{example}
When $p=3$, find a monic null polynomial of degree
$\omega_1(3^4)=3^2=9$ modulo $3^4$.
\end{example}
\begin{solution}
Choosing $q_3(x)=1$, the system of congruences becomes
\[
\left[\begin{matrix}%
1 & 1\\
2 & 2^2\\
\end{matrix}\right]
\left[\begin{matrix}%
-q_1(j)\\q_2(j)
\end{matrix}\right]\equiv
\left[\begin{matrix}%
1\\2\\
\end{matrix}\right]\pmod p\Rightarrow
\left[\begin{matrix}%
1 & 0\\
1 & 1\\
\end{matrix}\right]
\left[\begin{matrix}%
q_1(j)\\q_2(j)
\end{matrix}\right]\equiv
\left[\begin{matrix}%
2\\2\\
\end{matrix}\right]\pmod p.
\]
Solving this system of congruences, one has $q_1(j)\equiv 2\pmod p$
and $q_2(j)\equiv 0\pmod p$ for any integer $j$. So, choosing
$q_1(x)=2$ and $q_2(x)=0$, one gets a null polynomial
$f(x)=(x(x-1)(x-2))^3+3^2(x(x-1)(x-2))\cdot
2=\uwave{(x^3-3x^2+2x)^3+18(x^3-3x^2+2x)}\equiv\uwave{x^9-9x^8+33x^7+18x^6-15x^5-36x^4+26x^3+27x^2+36x}\pmod{3^4}$.
Experiments have been made to verify this null polynomial.
\end{solution}

Next, we give a least-degree monic null polynomial modulo $p^{p+1}$
by combining $\NPp$ directly, without solving the system of
congruences.

\begin{definition}
Assume $p$ is a prime. Define
$\NPpd{p+1}(x)=\prod_{i=0}^{p-1}(\NPp-ip)$. The following theorem
ensures that this polynomial is a least-degree monic null polynomial
modulo $p^{p+1}$, which plays an important role to construct
least-degree monic null polynomials modulo $p^d$, together with
$\NPp$, the null polynomial modulo $p$.
\end{definition}

\begin{theorem}\label{theorem:base-null-poly-p+1}
Assume $p$ is a prime and $d\geq 2$. Then, $\NPpd{p+1}(x)$ is a
least-degree monic null polynomial modulo $p^{p+1}$.
\end{theorem}
\begin{proof}
Consider
$\frac{\NPpd{p+1}(x)}{p^p}=\prod_{i=0}^{p-1}\left(\frac{\NPp}{p}-i\right)$.
Since $\frac{\NPp}{p}$ is an integer, there must exists
$i\in\{0,\cdots,p-1\}$ such that $\frac{\NPp}{p}-i\equiv 0\pmod p$.
So $\frac{\NPpd{p+1}(x)}{p^p}$ is a null polynomial modulo $p$ and
then $\NPpd{p+1}(x)$ is a null polynomial modulo $p^{p+1}$.
Considering that
$\degree(\NPpd{p+1}(x),p^{p+1})=\omega_1(p^{p+1})=p^2$, this theorem
is thus proved.
\end{proof}

\begin{lemma}
Assume $p$ is a prime. Then, $\forall i,x\in\mathbb{Z}$,
$\lfloor(ip^2+x)/p\rfloor\equiv\lfloor x/p\rfloor\pmod p$.
\end{lemma}
\begin{proof}
It is easy to prove this lemma as follows:
$\lfloor(ip^2+x)/p\rfloor=\lfloor ip+x/p\rfloor=ip+\lfloor
x/p\rfloor\equiv\lfloor x/p\rfloor\pmod p$.
\end{proof}
\begin{lemma}\label{lemma:floor-Lambda-level1}
Assume $p$ is a prime. Then, $\forall i,j,x\in\mathbb{Z}$, $\lfloor
j\Lambda(ip^2+x)/p\rfloor\equiv\lfloor j\Lambda(x)/p\rfloor\pmod p$.
\end{lemma}
\begin{proof}
One has $j\Lambda(ip^2+x)=j\prod_{0\leq k\leq p-1 \atop k\not\equiv
ip^2+x\pmod p}(ip^2+(x-k))=jA(ip^2+x)p^2+j\prod_{0\leq k\leq p-1
\atop k\not\equiv x\pmod p}(x-k)=jA(ip^2+x)p^2+j\Lambda(x)$, where
$A(ip^2+x)$ is the sum of all terms in $\Lambda(ip^2+x)$ that can be
divided by $p^2$. Then, from the above lemma, one immediately prove
$\lfloor j\Lambda(ip^2+x)/p\rfloor\equiv\lfloor
j\Lambda(x)/p\rfloor\pmod p$.
\end{proof}
\begin{theorem}\label{theorem:base-null-poly-p+1-feature}
Assume $p$ is a prime. Then, $\forall i\in\mathbb{Z}$ and $\forall
j_0,j_1\in\{0,\cdots,p-1\}$,
$\frac{\NPpd{p+1}(ip^2+j_1p+j_0)}{p^{p+1}}\equiv i-\lfloor
j_1\Lambda(j_1p+j_0)/p\rfloor\pmod p$.
\end{theorem}
\begin{proof}
Assuming $x=ip^2+j_1p+j_0$, similar to the proof of the above
theorem, one has
$\frac{\NPpd{p+1}(x)}{p^{p+1}}=\frac{\NPpf((ip+j_1)\Lambda(x))}{p}$.
Then, $\frac{\NPpd{p+1}(x)}{p^{p+1}}\equiv\left(i\Lambda(x)+\lfloor
j_1\Lambda(x)/p\rfloor\right)\Lambda((ip+j_1)\Lambda(x)))\equiv-(-i+\lfloor
j_1\Lambda(x)/p\rfloor)=i-\lfloor j_1\Lambda(x)/p\rfloor\equiv
i-\lfloor j_1\Lambda(j_1p+j_0)/p\rfloor\pmod p$. This immediately
proves this theorem.
\end{proof}
From the above theorem, one can easily derive the following
corollaries.
\begin{corollary}\label{corollary:base-null-poly-p+1-feature}
Assume $p$ is a prime. Then, $\forall j\in\{0,\cdots,p^2-1\}$,
$\left\{\frac{\NPpd{p+1}(ip^2+j)}{p^{p+1}}\right\}$ forms a complete
system of residues modulo $p$ when $i$ runs through a complete
system of residues modulo $p$.
\end{corollary}
\begin{corollary}\label{corollary:base-null-poly-p+1-feature2}
Assume $p$ is a prime. Then, $\forall i\in\mathbb{Z}$ and $\forall
j\in\{0,\cdots,p^3-1\}$,
$\frac{\NPpd{p+1}(ip^3+j)}{p^{p+1}}\equiv\frac{\NPpd{p+1}(j)}{p^{p+1}}\pmod
p$, i.e., $\frac{\NPpd{p+1}(x)}{p^{p+1}}$ forms a periodic function
modulo $p$ of period $p^3$.
\end{corollary}
\begin{corollary}\label{corollary:base-null-poly-p+1-feature3}
Assume $p$ is a prime. Then, $\forall i\in\mathbb{Z}$ and $\forall
j\in\{0,\cdots,p-1\}$, $\frac{\NPpd{p+1}(ip^2+j)}{p^{p+1}}\equiv
i\pmod p$.
\end{corollary}

\subsubsection[The Case of $p+2\leq d\leq 2p+1$]{The Case of $\bm{p+2\leq d\leq 2p+1}$}

\begin{lemma}\label{lemma:null-poly-p+2}
Assume $p$ is a prime. A polynomial $f(x)$ is a null polynomial
modulo $p^{p+2}$, if and only if
\[
f(x)\equiv \NPpd{p+1}(x)\NPp
q_{p+2}^*(x)+p\NPpd{p+1}(x)q_{p+1}^*(x)+\sum_{n=p}^1p^{p+2-n}(\NPp)^nq_n^*(x)\pmod{p^{p+2}},
\]
where $q_1^*(x)$, $\cdots$, $q_{p+1}^*(x)$ are any polynomials of
degree less than $p$ modulo $p^{p+2}$ and $q_{p+2}^*(x)$ is any
polynomial of any degree modulo $p^{p+2}$.
\end{lemma}
\begin{proof}
The ``if" part is obvious, so we only prove the ``only if" part.

From Corollary \ref{corollary:null-poly-hierarchy-pd}, one has
$f(x)\equiv
\NPpd{p+1}(x)q_{p+1}(x)+\sum_{n=p-1}^0p^{p+1-n}(\NPp)^nq_n(x)\pmod{p^{p+2}}$,
where $q_0(x)$, $\cdots$, $q_{p-1}(x)$ are polynomials of degree
less than $p$ modulo $p^{p+2}$.

Assuming $x=ip+j$, where $j\in\{0,\cdots,p-1\}$ and $i$ runs through
a complete system of residues modulo $p^{p+1}$. Substituting
$x=ip+j$ into $f(x)\equiv 0\pmod{p^{p+2}}$, one has
$\NPpd{p+1}(ip+j)q_{p+1}(ip+j)+\sum_{n=p-1}^0p^{p+1-n}(-ip)^nq_n(x)\equiv
0\pmod{p^{p+2}}\Rightarrow
\frac{\NPpd{p+1}(ip+j)}{p^{p+1}}q_{p+1}(j)+\sum_{n=p-1}^0(-i)^nq_n(j)\equiv
0\pmod p$. Choosing $i=0$, one has $\NPpd{p+1}(j)=0$ and
$\sum_{n=p-1}^1(-i)^nq_n(x)\equiv 0\pmod p$, so $q_0(j)\equiv 0\pmod
p$. Then, removing $q_0(j)$ from the congruence, one has
$\sum_{n=p-1}^1(-i)^nq_n(j)\equiv
-\frac{\NPpd{p+1}(ip+j)}{p^{p+1}}q_{p+1}(j)\pmod p$. From Theorem
\ref{theorem:base-null-poly-p+1-feature}, choosing $i_1\not\equiv
i_2\pmod p$, one has
$\frac{\NPpd{p+1}(i_1p^2+j)}{p^{p+1}}=i_1\not\equiv
i_2=\frac{\NPpd{p+1}(i_2p^2+j)}{p^{p+1}}\pmod p$. Then, one has
$(i_1-i_2)q_{p+1}(j)\equiv 0\pmod p$, which immediately leads to the
fact that $q_{p+1}(j)\equiv 0\pmod p$ and $q_1(j)\equiv\cdots\equiv
q_{p-1}(j)\equiv 0\pmod p$ hold for any integer
$j\in\{0,\cdots,p-1\}$. That is, all the $p$ polynomials are null
polynomials modulo $p$. In a similar way used in above proofs, we
can write each polynomial as two parts, and then prove this lemma.
\end{proof}

\begin{lemma}
Assume $p$ is a prime, then $\omega_1(p^{p+2})=p(p+1)$.
\end{lemma}
\begin{proof}
From the above lemma, choosing $q_{p+2}(x)=1$, this lemma is thus
proved.
\end{proof}

\begin{theorem}
Assume $p$ is a prime and $1\leq d\leq p$. A polynomial $f(x)$ is a
null polynomial modulo $p^{p+1+d}$, if and only if
\[
f(x)\equiv
\sum_{n_1=d}^0p^{d-n_1}\NPpd{p+1}(x)(\NPp)^{n_1}q_{p+n_1}^*(x)+\sum_{n_2=p}^1p^{p+1+d-n_2}(\NPp)^{n_2}q_{n_2}^*(x)\pmod{p^{p+1+d}},
\]
where $q_1(x)$, $\cdots$, $q_{p+d}(x)$ are any polynomials of degree
less than $p$ modulo $p^{p+1+d}$ and $q_{p+1+d}(x)$ is any
polynomial of any degree modulo $p^{p+1+d}$.
\end{theorem}
\begin{proof}
The ``if" part is obvious. We use induction on $d$ to prove the
``only if" part.

When $d=1$, this theorem has been proved in the above lemma. Under
the assumption that this theorem is true for any integer less than
$d$, let us prove the case of $d\geq 2$.

Note that this theorem means that $\forall 1\leq n\leq d-1$,
$\omega_1(p^{p+1+n})=p(p+n)$, so $\NPpd{p+1}(x)(\NPp)^n$ is a
least-degree monic null polynomial modulo $p^{p+1+n}$. Then, from
Corollary \ref{corollary:null-poly-hierarchy-pd}, one has
\[
f(x)\equiv
\sum_{n_1=d-1}^0p^{d-1-n_1}\NPpd{p+1}(x)(\NPp)^{n_1}q_{p+1+n_1}(x)+\sum_{n_2=p-1}^0p^{p+d-n_2}(\NPp)^{n_2}q_{n_2}(x)\pmod{p^{p+1+d}},
\]
where except $q_0(x)$, $\cdots$, $q_{p+d-1}(x)$ are of degree less
than $p$.

Assuming $x=ip+j$, where $j\in\{0,\cdots,p-1\}$ and $i$ runs through
a complete system of residues modulo $p^{p+1}$. Substituting
$x=ip+j$ into $f(x)\equiv 0\pmod{p^{p+2}}$, one has
\[
\sum_{n_1=d-1}^0p^{d-1-n_1}\NPpd{p+1}(ip+j)(-ip)^{n_1}q_{p+1+n_1}(ip+j)+\sum_{n_2=p-1}^0p^{p+d-n_2}(-ip)^{n_2}q_{n_2}(ip+j)\equiv
0\pmod{p^{p+1+d}},
\]
which leads to
$\frac{\NPpd{p+1}(ip+j)}{p^{p+1}}\sum_{n_1=d-1}^0(-i)^{n_1}q_{p+1+n_1}(j)+\sum_{n_2=p-1}^0(-i)^{n_2}q_{n_2}(j)\equiv
0\pmod p$. From Corollary
\ref{corollary:base-null-poly-p+1-feature}, $\forall
i,j\in\{0,\cdots,p-1\}$,
$\left\{\frac{\NPpd{p+1}(kp^2+ip+j)}{p^{p+1}}\right\}_{0\leq k\leq
p-1}$ forms a complete system of residues modulo $p$. This means
that
$\sum_{n_1=d-1}^0(-i)^{n_1}q_{p+1+n_1}(j)\equiv\sum_{n_2=p-1}^0(-i)^{n_2}q_{n_2}(j)\equiv
0\pmod p$. Since $1\leq d\leq p$, due to the same reason given in
the proof of Theorem \ref{theorem:null-poly-2-p}, one has
$q_0(j)\equiv\cdots\equiv q_{p+d}(j)\equiv 0\pmod p$. That is,
$q_0(x)$, $\cdots$, $q_{p+d}(x)$ are all null polynomials modulo
$p$. Noticing that $q_{p+d}(x)=\NPp q_{p+d+1}^*(x)+pq_{p+d}^*(x)$
and $q_n(x)=pq_n^*(x)$ when $0\leq n\leq p+d-1$, one can prove this
theorem.
\end{proof}

\begin{theorem}
Assume $p$ is a prime and $p+2\leq d\leq 2p+1$, then
$\omega_1(p^d)=p(d-1)$.
\end{theorem}
\begin{proof}
From the above theorem, choosing $q_d(x)=1$, this theorem is thus
proved.
\end{proof}

\subsubsection[The Case of $d=2p+2$]{The Case of $\bm{d=2p+2}$}

Apparently, the case of $d=2p+2$ is an analog of the case of
$d=p+1$. The following theorems can be easily obtained via the same
way as above, so the proofs are omitted here.

\begin{theorem}\label{theorem:null-poly-2p+2}
Assume $p$ is a prime and $1\leq d\leq p$. A polynomial $f(x)$ is a
null polynomial modulo $p^{2p+2}$, if and only if
\begin{itemize}
\item
$f(x)\equiv
\sum_{n_1=p}^0p^{p-n_1}\NPpd{p+1}(x)(\NPp)^{n_1}q_{p+1+n_1}(x)+\sum_{n_2=p}^1p^{2p+2-n_2}(\NPp)^{n_2}q_{n_2}^*(x)\pmod{p^{2p+2}}$;

\item
$q_1^*(x)$, $\cdots$, $q_p^*(x)$ are any polynomials of degree less
than $p$ modulo $p^{2p+2}$;

\item
$q_{2p+1}(x)$ is any polynomial of any degree and $q_{p+1}(x)$,
$\cdots$, $q_{2p}(x)$ polynomials of degree less than $p$ modulo
$p^{2p+2}$ that satisfy
$\sum_{n_1=p-1}^0(-i)^{n_1}q_{p+1+n_1}(j)\equiv iq_{2p+1}(j)\pmod p$
for $i=\lfloor x/p\rfloor$ and $j=x\bmod p$.
\end{itemize}
\end{theorem}

\begin{theorem}
Assume $p$ is a prime, then $\omega_1(p^{2p+2})=2p^2$.
\end{theorem}

\begin{corollary}
Assume $p$ is a prime. If $f(x)$ is a least-degree monic null
polynomial modulo $p^{p+1}$, then $(f(x))^2$ is a least-degree monic
null polynomial modulo $p^{2p+2}$.
\end{corollary}

\begin{remark}
Specially, $(\NPpd{p+1}(x))^2$ is a least-degree monic null
polynomial modulo $p^{2p+2}$.
\end{remark}

\subsubsection[The Case of $1\leq d\leq p(p+1)$]{The Case of $\bm{1\leq d\leq p(p+1)}$}

Generalizing the above procedure until $d=p(p+1)$, one has the
following theorems.

\begin{definition}
Assume $p$ is a prime. Define
$A_{p,i,s,t}(x)=\sum_{n_i=s}^tp^{s-n_i}(\NPp)^{n_i}q_{i,n_i}(x)$.
\end{definition}

\begin{theorem}
Assume $p$ is a prime and $1\leq k\leq p$. Then, $f(x)$ is a null
polynomial modulo $p^{k(p+1)}$ if and only if
\begin{multline*}
f(x)\equiv
(\NPpd{p+1}(x))^{k-1}A_{p,k-1,p,0}(x)+\sum_{i=k-2}^1p^{(k-1-i)(p+1)}(\NPpd{p+1}(x))^iA_{p,i,p-1,0}(x)\\
+p^{(k-1)(p+1)}A_{p,0,p-1,1}(x)\pmod{p^{k(p+1)}},
\end{multline*}
where all $q$-polynomials but the highest one in $A_{p,k-1,p,0}(x)$
is of degree less than $p$ modulo $p^{d(p+1)}$, and the
$q$-polynomials in $A_{p,k-1,p,0}(x)$ satisfy
$\sum_{n_{k-1}=p-1}^0(-i)^{n_{k-1}}q_{k-1,n_{k-1}}(j)\equiv
iq_{k-1,p}(j)\pmod p$ for $i=\lfloor x/p\rfloor$ and $j=x\bmod p$.
Specially, $(\NPpd{p+1}(x))^k$ is a least-degree monic null
polynomial modulo $p^{k(p+1)}$.
\end{theorem}

\begin{theorem}
Assume $p$ is a prime, $k(p+1)<d<(k+1)(p+1)$ and $d_k=d-k(p+1)$,
where $1\leq k\leq p$. Then, $f(x)$ is a null polynomial modulo
$p^d$ if and only if,
\begin{multline*}
f(x)\equiv
(\NPpd{p+1}(x))^{k-1}A_{p,k-1,d_k,0}(x)+\sum_{i=k-2}^1p^{(k-1-i)(p+1)}(\NPpd{p+1}(x))^iA_{p,i,p-1,0}(x)\\
+p^{(k-1)(p+1)}A_{p,0,p-1,1}(x)\pmod{p^{k(p+1)}},
\end{multline*}
where all $q$-polynomials but the highest one in $A_{p,k-1,p,0}(x)$
is of degree less than $p$ modulo $p^{d(p+1)}$, and all the
$q$-polynomials are null polynomials modulo $p$.
\end{theorem}

\begin{theorem}
Assume $p$ is a prime and $1\leq d\leq p(p+1)$, then
$\omega_1(p^d)=kp^2+p(d-k(p+1))=p(d-k)$, where
$k=\left\lfloor\frac{d}{p+1}\right\rfloor$.
\end{theorem}

The above theorems can be proved via a complicated induction on
$k\geq 1$: from $d=(k-1)(p-1)$ to $(k-1)(p-1)<d<k(p-1)$ and then to
$d=k(p-1)$, and finally to $k(p-1)<d<(k+1)(p-1)$. Since we will give
a similar proof on the general case of $d\geq 1$ later, the proofs
of the above theorems are omitted here.

\subsubsection[The Case of $d=p(p+1)+1$]{The Case of $\bm{d=p(p+1)+1}$}

\begin{theorem}
Assume $p$ is a prime. Then, $f(x)$ is a null polynomial modulo
$p^{p(p+1)+1}$ if and only if
\begin{multline*}
f(x)\equiv
(\NPpd{p+1}(x))^pQ(x)+\sum_{i=p-1}^1p^{(k-1-i)(p+1)}(\NPpd{p+1}(x))^iA_{p,i,p-1,0}(x)\\
+p^{(k-1)(p+1)}A_{p,0,p-1,1}(x)\pmod{p^{k(p+1)}},
\end{multline*}
where all $q$-polynomials are of degree less than $p$ modulo
$p^{d(p+1)}$, $Q(x)\equiv 1\pmod p$, the $q$-polynomials in
$A_{p,1,p,0}(x)$ satisfy
$\sum_{n_1=p-1}^0(-i)^{n_1}q_{k-1,n_1}(j)\equiv 1\pmod p$ for
$i=\lfloor x/p\rfloor$ and $j=x\bmod p$, and all other
$q$-polynomials are null polynomials modulo $p$.
\end{theorem}

\begin{theorem}
Assume $p$ is a prime, then
$\omega_1(p^{p(p+1)})=\omega_1(p^{p(p+1)+1})=p^3$.
\end{theorem}

\begin{definition}
Using the same way of defining $\NPpd{p+1}(x)$, we have
\[
\NPpd{p(p+1)+1}(x)\equiv
\prod_{i=0}^{p-1}\left(\NPpd{p+1}(x)-ip^{p+1}\right)\equiv
0\pmod{p^{p(p+1)+1}},
\]
i.e.,
\[
\frac{\NPpd{p(p+1)+1}(x)}{p^{p(p+1)}}=\prod_{i=0}^{p-1}\left(\frac{\NPpd{p+1}(x)}{p^{p+1}}-i\right)\equiv
0\pmod p.
\]
\end{definition}

\begin{theorem}
Assume $p$ is a prime. Then, $\NPpd{p(p+1)+1}(x)$ is a least-degree
monic null polynomial modulo $p^{p(p+1)+1}$.
\end{theorem}
\begin{proof}
The same as the proof of Theorem \ref{theorem:base-null-poly-p+1}.
\end{proof}

\begin{lemma}
Assume $p$ is a prime. Then, $\forall i_1,i_2,j_1,j_2\in\mathbb{Z}$,
$\lfloor\lfloor
(i_1p^3+j_1)/p\rfloor(i_2p^2+j_2)/p\rfloor\equiv\lfloor\lfloor
j_1/p\rfloor j_2/p\rfloor\pmod p$.
\end{lemma}
\begin{proof}
One has $\lfloor\lfloor
(i_1p^3+j_1)/p\rfloor(i_2p^2+j_2)/p\rfloor=\lfloor\lfloor
i_1p^2+j_1/p\rfloor(i_2p+j_2/p)\rfloor=\lfloor(i_1p^2+\lfloor
j_1/p\rfloor)(i_2p+j_2/p)\rfloor=\lfloor
i_1i_2p^3+i_1j_2p+i_2p\lfloor j_1/p\rfloor+\lfloor j_1/p\rfloor
j_2/p\rfloor=i_1i_2p^3+i_1j_2p+i_2p\lfloor
j_1/p\rfloor+\lfloor\lfloor j_1/p\rfloor
j_2/p\rfloor\equiv\lfloor\lfloor j_1/p\rfloor j_2/p\rfloor\pmod p$.
\end{proof}
\begin{lemma}
Assume $p$ is a prime. Then, $\forall i_1,i_2,j_1,j_2\in\mathbb{Z}$,
$\lfloor\lfloor
\Lambda(i_1p^3+j_1)/p\rfloor\Lambda(i_2p^2+j_2)/p\rfloor\equiv\lfloor\lfloor
\Lambda(j_1)/p\rfloor\Lambda(j_2)/p\rfloor\pmod p$ and
$\lfloor\Lambda(i_1p^2+j_1)/p\Lambda(i_2p^2+j_2)/p\rfloor\equiv
\lfloor\Lambda(j_1)/p\Lambda(j_2)/p\rfloor\pmod p$.
\end{lemma}
\begin{proof}
This lemma can be proved in a similar to Lemma
\ref{lemma:floor-Lambda-level1}, based on the above lemma.
\end{proof}
\begin{theorem}\label{theorem:base-null-poly-p+1-feature}
Assume $p$ is a prime. Then, $\forall i\in\mathbb{Z}$ and $\forall
j_0,j_1,j_2\in\{0,\cdots,p-1\}$,
$\frac{\NPpd{p(p+1)+1}(ip^3+j_2p^2+j_1p+j_0)}{p^{p(p+1)+1}}\equiv-(i+\lfloor
j_2\Lambda(j_{1,0})\Lambda(j_{2,1})/p\rfloor+\lfloor\lfloor
j_1\Lambda(j_{1,0})/p\rfloor\Lambda(j_{2,1})/p\rfloor)\pmod p$,
where $j_{1,0}=j_1p+j_0$ and $j_{2,1}=(j_2p+j_1)\Lambda(j_1p+j_0)$.
\end{theorem}
\begin{proof}
Assuming $x_0=ip^3+j_2p^2+j_1p+j_0$, one has
$\frac{\NPpd{p+1}(x)}{p^{p+1}}=\frac{\NPpf((ip^2+j_2p+j_1)\Lambda(x_0))}{p}$.
Then, assuming $x_1=(ip^2+j_2p+j_1)\Lambda(x_0)$, one has
$\frac{\NPpd{p+1}(x)}{p^{p+1}}=((ip+j_2)\Lambda(x_0)+\lfloor
j_1\Lambda(x_0)/p\rfloor)\Lambda(x_1)=x_2$. Next,
$\frac{\NPpd{p(p+1)+1}(ip^3+j_2p^2+j_1p+j_0)}{p^{p(p+1)+1}}=\left(i\Lambda(x_0)\Lambda(x_1)+\lfloor
j_2\Lambda(x_0)\Lambda(x_1)/p\rfloor+\lfloor\lfloor
j_1\Lambda(x_0)/p\rfloor\Lambda(x_1)/p\rfloor\right)\Lambda(x_2)$.
Then, from the above lemma, one has
$\frac{\NPpd{p+1}(x)}{p^{p+1}}\equiv-(i+\lfloor
j_2\Lambda(x_0)\Lambda(x_1)/p\rfloor+\lfloor\lfloor
j_1\Lambda(x_0)/p\rfloor\Lambda(x_1)/p\rfloor)\equiv-(i+\lfloor
j_2\Lambda(j_{1,0})\Lambda(j_{2,1})/p\rfloor+\lfloor\lfloor
j_1\Lambda(j_{1,0})/p\rfloor\Lambda(j_{2,1})/p\rfloor)\pmod p$. Thus
this theorem is proved.
\end{proof}
From the above theorem, one can easily derive the following
corollaries.
\begin{corollary}
Assume $p$ is a prime. Then, $\forall j\in\{0,\cdots,p^3-1\}$,
$\left\{\frac{\NPpd{p(p+1)+1}(ip^3+j)}{p^{p(p+1)+1}}\right\}$ forms
a complete system of residues modulo $p$ when $i$ runs through a
complete system of residues modulo $p$.
\end{corollary}
\begin{corollary}
Assume $p$ is a prime. Then, $\forall i\in\mathbb{Z}$ and $\forall
j\in\{0,\cdots,p^4-1\}$,
$\frac{\NPpd{p(p+1)+1}(ip^4+j)}{p^{p(p+1)+1}}\equiv\frac{\NPpd{p+1}(j)}{p^{p+1}}\pmod
p$, i.e., $\frac{\NPpd{p+1}(x)}{p^{p+1}}$ forms a function modulo
$p$ of period $p^4$.
\end{corollary}
\begin{corollary}
Assume $p$ is a prime. Then, $\forall i\in\mathbb{Z}$ and $\forall
j\in\{0,\cdots,p-1\}$,
$\frac{\NPpd{p(p+1)+1}(ip^3+j)}{p^{p(p+1)+1}}\equiv -i\pmod p$.
\end{corollary}

Observing $\NPp$, $\NPpd{p+1}(x)$ and $\NPpd{p(p+1)+1}(x)$ and
comparing their features, one can find a recursive formula to
further generalize the above results to all values of $d$. This
leads to an inductive proof of the general case of $d\geq 1$ given
in next subsection.

\subsection[The General Case: Null Polynomials modulo $p^d$ ($d\geq 1$)]%
{The General Case: Null Polynomials modulo $\bm{p^d}$ ($\bm{d\geq
1}$)} \label{subsection:null-poly-pd}

\subsubsection{Some Definitions and Preliminaries}

\begin{definition}
Assume $p$ is a prime and $n\geq 0$. Define an index-sequence as
follows:
\[
I_p(n)=\begin{cases}%
0, & \mbox{when }n=0,\\
pI_p(n-1)+1, & \mbox{when }n\geq 1.
\end{cases}
\]
When $n\geq 1$, one can easily derive the close form of the above
sequence: $I_p(n)=\sum_{i=0}^{n-1}p^i=\frac{p^n-1}{p-1}$. Specially,
$I_p(1)=1$, $I_p(2)=p+1$, $I_p(3)=p(p+1)+1=p^2+p+1$.
\end{definition}
\begin{definition}
Assume $p$ is a prime and $n\geq 0$. Define a monic integer
polynomial as follows:
\[
\mathcal{G}_{p,n}(x)=\begin{cases}%
x, & \mbox{when }n=0\\
\prod_{i=0}^{p-1}\left(\mathcal{G}_{p,n-1}(x)-ip^{I_p(n-1)}\right),
& \mbox{when }n\geq 1,
\end{cases}
\]
where $I_p(n)$ denotes the above index-sequence. Specially,
$\mathcal{G}_{p,1}(x)=\NPp$,
$\mathcal{G}_{p,2}(x)=\NPpd{p+1}(x)=\NPpd{I_p(2)}(x)$ and
$\mathcal{G}_{p,3}(x)=\NPpd{p(p+1)+1}(x)=\NPpd{I_p(3)}(x)$.
\end{definition}

An equivalent definition of the above polynomial is as follows.

\begin{definition}
Assume $p$ is a prime and $n\geq 0$. Define a rational
polynomial\footnote{Like the name of ``integer polynomial", a
``rational polynomial" means a polynomial with rational
coefficients.} $\widetilde{\mathcal{G}}_{p,n}(x)$ as follows:
\[
\widetilde{\mathcal{G}}_{p,n}(x)=\begin{cases}%
x, & \mbox{when }n=0\\
\frac{\NPpf\left(\widetilde{\mathcal{G}}_{p,n-1}(x)\right)}{p}=
\frac{\prod_{i=0}^{p-1}\left(\widetilde{\mathcal{G}}_{p,n-1}(x)-i\right)}{p},
& \mbox{when }n\geq 1.
\end{cases}
\]
Then, define an integer polynomial by
$\mathcal{G}_{p,n}(x)=p^{I_p(n)}\widetilde{\mathcal{G}}_{p,n}(x)$.
\end{definition}

In the following, both $\mathcal{G}_{p,n}(x)$ and
$\widetilde{\mathcal{G}}_{p,n}(x)$ will be frequently used to
achieve a more concise description of the results on polynomials of
this kind.

\begin{theorem}\label{theorem:base-null-poly-feature}
Assume $p$ is a prime and $n\geq 1$. Then, $\forall
x=ip^n+\sum_{k=n-1}^0j_ip^i$, where $i\in\mathbb{Z}$ and
$j_0,\cdots,j_{n-1}\in\{0,\cdots,p-1\}$, it is true that
$\widetilde{\mathcal{G}}_{p,n}(x)=\frac{\mathcal{G}_{p,n}(x)}{p^{I_p(n)}}=
i\left(\prod_{k=0}^{n-1}\Lambda\left(\widetilde{\mathcal{G}}_{p,k}(x)\right)+A_n(x)p\right)+J_n(j_0,\cdots,j_{n-1})$.
\end{theorem}
\begin{proof}
It is obvious that the result is true when $n=1$. Let us use
induction on $n$ to prove the general case of $n\geq 2$, under the
assumption that this theorem is true for any integer less than $n$.

From the hypothesis,
$\widetilde{\mathcal{G}}_{p,n-1}(x)=(ip+j_{n-1})\left(\prod_{k=0}^{n-2}\Lambda\left(\widetilde{\mathcal{G}}_{p,k}(x)\right)+A_{n-1}(x)p\right)+J_{n-1}(j_0,\cdots,j_{n-2})$.
Then, from the definition of $\widetilde{\mathcal{G}}_{p,n}(x)$ and
Theorem \ref{theorem:base-null-poly-p-form2}, one has
\begin{eqnarray*}
\widetilde{\mathcal{G}}_{p,n}(x) & = &
\frac{\NPpf\left(\widetilde{\mathcal{G}}_{p,n-1}(x)\right)}{p}
=\left\lfloor\widetilde{\mathcal{G}}_{p,n-1}(x)/p\right\rfloor\Lambda\left(\widetilde{\mathcal{G}}_{p,n-1}(x)\right)\\
& = &
\left\lfloor\left.\left((ip+j_{n-1})\left(\prod_{k=0}^{n-2}\Lambda\left(\widetilde{\mathcal{G}}_{p,k}(x)\right)+A_{n-1}(x)p\right)+J_{n-1}(j_0,\cdots,j_{n-2})\right)\right/p\right\rfloor\Lambda\left(\widetilde{\mathcal{G}}_{p,n-1}(x)\right)\\
& = & i\left(\prod_{k=0}^{n-1}\Lambda\left(\widetilde{\mathcal{G}}_{p,k}(x)\right)+A_{n-1}\Lambda\left(\widetilde{\mathcal{G}}_{p,n-1}(x)\right)p\right)\\
& & +\left(A_{n-1}j_{n-1}+\left\lfloor\left.\left(
j_{n-1}\prod_{k=0}^{n-2}\Lambda\left(\widetilde{\mathcal{G}}_{p,k}(x)\right)+J_{n-1}(j_0,\cdots,j_{n-2})\right)\right/p\right\rfloor\right)\Lambda\left(\widetilde{\mathcal{G}}_{p,n-1}(x)\right).
\end{eqnarray*}
Applying the hypothesis on $\widetilde{\mathcal{G}}_{p,k}(x)$, one
has
$\widetilde{\mathcal{G}}_{p,k}(x)=(ip^{n-k}+\cdots+j_k)\left(\prod_{l=0}^{k-1}\Lambda\left(\widetilde{\mathcal{G}}_{p,l}(x)\right)+A_k(x)p\right)+J_k(j_0,\cdots,j_{k-1})$.
So, when $0\leq k\leq n-2$, we can get
$\widetilde{\mathcal{G}}_{p,k}(x)=ip^2D_k(x)+J_k^*(j_0,\cdots,j_{n-1})$
and then
$\Lambda\left(\widetilde{\mathcal{G}}_{p,k}(x)\right)=ip^2D_k^*(x)+\Lambda(J_k^*(j_0,\cdots,j_{n-1}))$.
Thus, moving all terms of $p^2$ and of high powers of $p$ out of the
floor function, one has
\[
\left\lfloor\left.\left(
j_{n-1}\prod_{k=0}^{n-2}\Lambda\left(\widetilde{\mathcal{G}}_{p,k}(x)\right)+J_{n-1}(j_0,\cdots,j_{n-2})\right)\right/p\right\rfloor
=ipD^*(x)+J_n^*(j_0,\cdots,j_{n-1}).
\]
Substituting the above equation into
$\widetilde{\mathcal{G}}_{p,n}(x)$, one immediately has
\begin{eqnarray*}
\widetilde{\mathcal{G}}_{p,n}(x) & = & i\left(\prod_{k=0}^{n-1}\left(\widetilde{\mathcal{G}}_{p,k}(x)\right)+A_{n-1}\Lambda\left(\widetilde{\mathcal{G}}_{p,n-1}(x)\right)p\right)\\
& &
+\left(A_{n-1}j_{n-1}+ipD^*(x)+J_n^*(j_0,\cdots,j_{n-1})\right)\Lambda\left(\widetilde{\mathcal{G}}_{p,n-1}(x)\right)\\
& = & i\left(\prod_{k=0}^{n-1}\Lambda\left(\widetilde{\mathcal{G}}_{p,k}(x)\right)+(A_{n-1}+D^*(x))\Lambda\left(\widetilde{\mathcal{G}}_{p,n-1}(x)\right)p\right)\\
& &
{}+\left(A_{n-1}j_{n-1}+J_n^*(j_0,\cdots,j_{n-1})\right)\Lambda\left(\widetilde{\mathcal{G}}_{p,n-1}(x)\right).
\end{eqnarray*}
Similarly, one has
$\Lambda\left(\widetilde{\mathcal{G}}_{p,n-1}(x)\right)=ipD_{n-1}^*(x)+\Lambda(J_{n-1}^*(j_0,\cdots,j_{n-1}))$.
Then,
\begin{eqnarray*}
g_n(x)& = &
i\left(\prod_{k=0}^{n-1}\Lambda\left(\widetilde{\mathcal{G}}_{p,k}(x)\right)+(A_{n-1}+D^*(x))\Lambda\left(\widetilde{\mathcal{G}}_{p,n-1}(x)\right)p\right)\\
& &
{}+\left(A_{n-1}j_{n-1}+J_n^*(j_0,\cdots,j_{n-1})\right)(ipD_{n-1}^*(x)+\Lambda(J_{n-1}^*(j_0,\cdots,j_{n-1})))\\
& = &
i\left(\prod_{k=0}^{n-1}\Lambda\left(\widetilde{\mathcal{G}}_{p,k}(x)\right)+\left((A_{n-1}+D^*(x))\Lambda\left(\widetilde{\mathcal{G}}_{p,n-1}(x)\right)+(A_{n-1}j_{n-1}+J_n^*(j_0,\cdots,j_{n-1}))D_{n-1}^*(x)\right)p\right)\\
& &
{}+(A_{n-1}j_{n-1}+J_n^*(j_0,\cdots,j_{n-1}))\Lambda(J_{n-1}^*(j_0,\cdots,j_{n-1})).
\end{eqnarray*}
Assigning
$A_n(x)=((A_{n-1}+D^*(x))\Lambda(g_{n-1}(x))+(A_{n-1}j_{n-1}+J_n^*(j_0,\cdots,j_{n-1}))D_{n-1}^*(x))$
and
$J_n(j_0,\cdots,j_{n-1})=(A_{n-1}j_{n-1}+J_n^*(j_0,\cdots,j_{n-1}))\Lambda(J_{n-1}^*(j_0,\cdots,j_{n-1}))$,
one immediately gets
\[
\widetilde{\mathcal{G}}_{p,n}(x)=i\left(\prod_{k=0}^{n-1}\Lambda\left(\widetilde{\mathcal{G}}_{p,k}(x)\right)+A_n(x)p\right)+J_n(j_0,\cdots,j_{n-1}).
\]
Thus, this theorem is proved.
\end{proof}
The above theorem immediately derives the following corollaries.
\begin{corollary}\label{corollary:base-null-polys-feature}
Assume $p$ is a prime and $n\geq 1$. Then, $\forall i\in\mathbb{Z}$
and $j\in\{0,\cdots,p^n-1\}$, it is true that
$\widetilde{\mathcal{G}}_{p,n}(ip^n+j)\equiv(-1)^ni+J_n(j)\pmod p$.
\end{corollary}
\begin{corollary}\label{corollary:base-null-polys-feature1}
Assume $p$ is a prime and $n\geq 1$. Then, $\forall
j\in\{0,\cdots,p^n-1\}$,
$\left\{\widetilde{\mathcal{G}}_{p,n}(ip^n+j)\right\}$ forms a
complete system of residues modulo $p$ when $i$ runs through a
complete system of residues modulo $p$.
\end{corollary}
\begin{corollary}\label{corollary:base-null-polys-feature2}
Assume $p$ is a prime and $n\geq 1$. Then, $\forall i\in\mathbb{Z}$
and $\forall j\in\{0,\cdots,p^{n+1}-1\}$,
$\widetilde{\mathcal{G}}_{p,n}(ip^{n+1}+j)\equiv
\widetilde{\mathcal{G}}_{p,n}(j)\pmod p$, i.e.,
$\widetilde{\mathcal{G}}_{p,n}(x)$ forms a periodic function modulo
$p$ of period $p^{n+1}$.
\end{corollary}

Specially, when $j_1=\cdots=j_{n-1}=0$, we have a much simpler (but
not so useful as one can see later) form of the above results.

\begin{theorem}
Assume $p$ is a prime and $n\geq 1$. Then, $\forall x=ip^n+j$, where
$i\in\mathbb{Z}$ and $j\in\{0,\cdots,p-1\}$, it is true that
$\widetilde{\mathcal{G}}_{p,n}(x)=i\prod_{k=0}^{n-1}\Lambda\left(\widetilde{\mathcal{G}}_{p,k}(x)\right)$.
\end{theorem}
\begin{proof}
Let us use induction on $n$ to prove this lemma. When $n=1$, one can
easily get
$\widetilde{\mathcal{G}}_{p,1}(x)=\frac{\NPpf(ip+j)}{p}=i\Lambda(x)$.
Then, assuming this lemma is true for any integer less than $n$, let
us prove the case of $n\geq 2$. From the definition of
$\widetilde{\mathcal{G}}_{p,n}(x)$, one has
$\widetilde{\mathcal{G}}_{p,n}(x)=\frac{\NPpf\left(\widetilde{\mathcal{G}}_{p,n-1}(ip^n+j)\right)}{p}
=\frac{\NPpf\left(ip\prod_{j=1}^{n-1}\Lambda(c_j)\right)}{p}=\left\lfloor\left.
ip\prod_{k=0}^{n-2}\Lambda\left(\widetilde{\mathcal{G}}_{p,k}(x)\right)\right/p\right\rfloor\Lambda\left(\widetilde{\mathcal{G}}_{p,n-1}(ip^n+j)\right)
=i\prod_{k=0}^{n-1}\Lambda\left(\widetilde{\mathcal{G}}_{p,k}(x)\right)$.
Thus, this lemma is proved.
\end{proof}
\begin{corollary}\label{corollary:base-null-poly-feature1}
Assume $p$ is a prime and $n\geq 1$. Then, $\forall i\in\mathbb{Z}$
and $j\in\{0,\cdots,p-1\}$,
$\widetilde{\mathcal{G}}_{p,n}\equiv(-1)^ni\pmod p$.
\end{corollary}

\begin{lemma}\label{lemma:G-is-null-polys}
Assume $p$ is a prime and $n\geq 1$. Then, $\mathcal{G}_{p,n}(x)$ is
a monic null polynomial of degree $p^n$ modulo $p^{I_p(n)}$, but not
a null polynomial modulo $p^{I_p(n)+1}$.
\end{lemma}
\begin{proof}
The second part of this lemma is a straightforward result of
Corollary \ref{corollary:base-null-poly-feature1}, since
$\mathcal{G}_{p,n}(ip^n+j)\equiv (-1)^nip^{I_p(n)}\not\equiv
0\pmod{p^{I_p(n)+1}}$ when $i\not\equiv 0\pmod p$.

The first part of this lemma can be proved via induction on $n$.
When $n=1$, the result is obviously true. Let us consider the case
of $n\geq 2$. Following the definition of $\mathcal{G}_{p,n}(x)$,
one has
$\frac{\mathcal{G}_{p,n}(x)}{p^{I_p(n)-1}}=\frac{\mathcal{G}_{p,n}(x)}{p^{pI_p(n-1)}}
=\prod_{i=0}^{p-1}\left(\frac{\mathcal{G}_{p,n-1}(x)}{p^{I_p(n-1)}}-i\right)$.
From the hypothesis, $\mathcal{G}_{p,n-1}(x)$ is a null polynomial
modulo $p^{I_p(n-1)}$, so
$\frac{\mathcal{G}_{p,n-1}(x)}{p^{I_p(n-1)}}\in\mathbb{Z}$ and then
$\frac{\mathcal{G}_{p,n}(x)}{p^{I_p(n)-1}}$ is a null polynomial
modulo $p$. This means $\mathcal{G}_{p,n}(x)$ is a null polynomial
modulo $p^{I_p(n)}$. In addition, it is obvious that
$\degree(\mathcal{G}_{p,n},p^{I_p(n)})=p\cdot\degree(\mathcal{G}_{p,n-1},p^{I_p(n-1)})=pp^{n-1}=p^n$.
Thus this lemma is proved.
\end{proof}
\begin{corollary}
Assume $p$ is a prime, $e_1$, $\cdots$, $e_n\geq 1$ and
$d=\sum_{i=1}^ne_iI(i)$. Then,
$f(x)=\prod_{i=1}^n(\mathcal{G}_{p,i}(x))^{e_i}$ is a monic null
polynomial modulo $p^d$.
\end{corollary}
\begin{proof}
This corollary can be easily derived from the above lemma.
\end{proof}

\begin{definition}
Assume $p$ is a prime and $d\geq 0$. Define a polynomial sequence
$\left\{\mathcal{H}_{p,d}(x)=\prod_{i=1}^{\infty}(\mathcal{G}_{p,i}(x))^{e_{d,i}}\right\}_{d\geq
0}$, where $0\leq e_{d,i}\leq p$, in the following recursive way:
\[
\mathcal{H}_{p,d}(x)=\begin{cases}%
1, & \mbox{when }d=0,\\
\mathcal{H}_{p,d-1}(x)\mathcal{G}_{p,1}(x), & \mbox{when }d\geq
1\mbox{ and }\max_{i=1}^{\infty}(e_{d-1,i})\leq p-1,\\
\dfrac{\mathcal{H}_{p,d-1}(x)\mathcal{G}_{p,i+1}(x)}{(\mathcal{G}_{p,i}(x))^p},
& \mbox{when }d\geq 1\mbox{ and }\exists i\mbox{ such that
}e_{d-1,i}=p,e_{d-1,1}=\cdots=e_{d-1,i-1}=0.
\end{cases}
\]
One can easily verify that in each polynomial defined in the above
rule, there exists at most one exponent that satisfies $e_{d,i}=p$
and all exponents after this exponent are zeros. Denote
$\max_{i=1}^{\infty}(e_{d,i})$, i.e., the maximal exponent of
$\mathcal{H}_{p,d}(x)$, by $E_{\max}(\mathcal{H}_{p,d})$.
\end{definition}

Apparently, $\mathcal{H}_{p,d}(x)$ is an analog of an integer
represented with radix $p+1$, except that the last non-zero exponent
of the former may be $p$. So we define a numeric representation of
$\mathcal{H}_{p,d}(x)$ with floating radix.

\begin{definition}
\uline{The numeric representation} of a polynomial
$\mathcal{H}_{p,d}(x)=\prod_{i=1}^{\infty}(\mathcal{G}_{p,i}(x))^{e_{d,i}}$
is an integer defined by
\[
\langle\mathcal{H}_{p,d}\rangle=\sum_{i=1}^{\infty}e_{d,i}I_p(i)=(\cdots,e_{d,1})_{I_p},
\]
where the subscript ``$I_p$" denotes the floating radix of
$\langle\mathcal{H}_{p,d}\rangle$ as an integer. We call $e_{d,i}$
\uline{the $i$-th digit of the polynomial}. If
$\langle\mathcal{H}_{p,d}\rangle=(e_{i_2},\cdots,e_{i_1})_{I_p}=\sum_{i=i_1}^{i_2}e_{d,i}I_p(i)$,
where $1\leq i_1\leq i_2$, we say $i_2$ is the \uline{digit length}
and $e_{i_1},e_{i_2}$ are the \uline{LSD (least significant digit)}
and \uline{MSD (most significant digit)} of the polynomial,
respectively.
\end{definition}

\iffalse
\begin{lemma}
Assume $p$ is a prime and $d\geq 1$, then
$\langle\mathcal{H}_{p,d+1}\rangle>\langle\mathcal{H}_{p,d}\rangle$.
\end{lemma}
\begin{proof}
It is obvious that the digit length of $\mathcal{H}_{p,d}(x)$ is
finite, since $d$ is finite. Assume
$\mathcal{H}_{p,d}(x)=\prod_{i=1}^n(\mathcal{G}_{p,i}(x))^{d_i}$.
Then, consider two different cases. When
$E_{\max}(\mathcal{H}_{p,d})\leq p-1$, it is obvious that
$\langle\mathcal{H}_{p,d+1}\rangle=\langle\mathcal{H}_{p,d}\rangle+1>\langle\mathcal{H}_{p,d}\rangle$.
When $E_{\max}(\mathcal{H}_{p,d})=p$, i.e., $\exists i$ such that
$e_{d,i}=p$ and $e_{d,1}=\cdots=e_{d,i-1}=0$, one has
$\langle\mathcal{H}_{p,d+1}\rangle=\langle\mathcal{H}_{p,d}\rangle-p\cdot
(p+1)^{i-1}+(p+1)^i=\langle\mathcal{H}_{p,d}\rangle+(p+1)^{i-1}>\langle\mathcal{H}_{p,d}\rangle$
(note that $i\geq 1$). Thus this lemma is proved.
\end{proof}
\begin{corollary}
Assume $p$ is a prime and $d_1>d_2\geq 1$, then
$\langle\mathcal{H}_{p,d_1}\rangle>\langle\mathcal{H}_{p,d_2}\rangle$.
\end{corollary}
\fi

\begin{lemma}\label{lemma:monic-null-poly-degrees}
Assume $p$ is a prime and $d\geq 1$. Then,
\begin{enumerate}
\item
when $E_{\max}(\mathcal{H}_{p,d})\leq p-1$,
$\degree(\mathcal{H}_{p,d+1},p)=\degree(\mathcal{H}_{p,d},p)+p$;

\item
when $E_{\max}(\mathcal{H}_{p,d})=p$,
$\degree(\mathcal{H}_{p,d+1},p)=\degree(\mathcal{H}_{p,d},p)$.
\end{enumerate}
\end{lemma}
\begin{proof}
It is obvious that the digit length of $\mathcal{H}_{p,d}(x)$ is
finite, since $d$ is finite. Assume
$\mathcal{H}_{p,d}(x)=\prod_{i=1}^n(\mathcal{G}_{p,i}(x))^{d_i}$.
Then, consider two different cases. When
$E_{\max}(\mathcal{H}_{p,d})\leq p-1$,
$\degree(\mathcal{H}_{p,d+1},p)=\degree(\mathcal{H}_{p,d},p)+p>\degree(\mathcal{H}_{p,d},p)$.
When $E_{\max}(\mathcal{H}_{p,d})=p$, i.e., $\exists i$ such that
$e_{d,i}=p$ and $e_{d,1}=\cdots=e_{d,i-1}=0$,
$\degree(\mathcal{H}_{p,d+1},p)=\degree(\mathcal{H}_{p,d},p)$ since
$\degree(\mathcal{G}_{p,i+1}(x),p)=\degree((\mathcal{G}_{p,i}(x))^p,p)=pI(i)$.
Thus this lemma is proved.
\end{proof}

\begin{lemma}
Assume $p$ is a prime and $d\geq 1$. Then, $\mathcal{H}_{p,d}(x)$ is
a monic null polynomial modulo $p^d$.
\end{lemma}
\begin{proof}
Let us use induction on $d$ to prove this lemma. When $n=1$, it is
obviously true. Assuming it is also true for any integer less than
$d$, consider the case of $d\geq 2$.

When $E_{\max}(\mathcal{H}_{p,d})\leq p-1$, one can see
$\mathcal{H}_{p,d}(x)=\mathcal{H}_{p,d-1}(x)\mathcal{G}_{p,1}(x)$ is
a null polynomial modulo $p^{d-1+1}=p^d$. When
$E_{\max}(\mathcal{H}_{p,d})=p$, i.e., $\exists i$ such that
$e_{d-1,i}=p$ and $e_{d-1,1}=\cdots=e_{d-1,i-1}=0$,
$\frac{\mathcal{H}_{p,d-1}(x)}{(\mathcal{G}_{p,i}(x))^p}$ is a null
polynomial modulo $p^{(d-1)-pI(i)}$ and $\mathcal{G}_{p,i+1}(x)$ is
a null polynomial modulo $p^{I(i+1)}=p^{pI(i)+1}$. Thus,
$\mathcal{H}_{p,d}(x)=\frac{\mathcal{H}_{p,d-1}(x)\mathcal{G}_{p,i+1}(x)}{(\mathcal{G}_{p,i}(x))^p}$
is a null polynomial modulo $p^{(d-1)-pI(i)+pI(i)+1}=p^d$ and this
lemma is proved.
\end{proof}
\begin{theorem}
Assume $p$ is a prime and $d\geq 1$. Then,
$\langle\mathcal{H}_{p,d}\rangle=d$.
\end{theorem}
\begin{proof}
This theorem can be proved in the same way as the above theorem.
\end{proof}

\subsubsection{Main Results}
\label{subsubsection:MainResults}

\begin{theorem}\label{theorem:base-null-poly-hierarchy}
Assume $p$ is a prime, $n\geq 1$. If $\sum_{0\leq e_i\leq p-1 \atop
0\leq i\leq
n}\left(\prod_{i=1}^n\left(\widetilde{\mathcal{G}}_{p,i}(x)\right)^{e_i}q_{e_1,\cdots,e_n}(x)\right)\equiv
0\pmod p$ holds for any integer $x\in\mathbb{Z}$, then $\forall
e_1,\cdots,e_n\in\{0,\cdots,p-1\}$, $q_{e_1,\cdots,e_n}(x)$ is a
null polynomial modulo $p$.
\end{theorem}
\begin{proof}
Assuming $x=ip^n+\sum_{k=n-1}^0j_kp^k\in\mathbb{Z}$, where
$i\in\mathbb{Z}$ and $j_0,\cdots,j_{n-1}\in\{0,\cdots,p-1\}$, let us
prove this theorem via induction on $n$.

When $n=1$, one has $x=ip+j_0$ and
$\widetilde{\mathcal{G}}_{p,1}(x)\equiv -i\pmod p$, and then
$\sum_{0\leq e_1\leq
p-1}\left(\widetilde{\mathcal{G}}_{p,1}(x)\right)^{e_1}q_{e_1}(x)\equiv\sum_{0\leq
e_1\leq p-1}(-i)^{e_1}q_{e_1}(j_0)\equiv 0\pmod p$. Consider the
polynomial $f(y)=\sum_{0\leq e_1\leq p-1}q_{e_1}(j_0)y^{e_1}$, one
can see that $f(y)\equiv 0\pmod p$ for any integer
$y=-i\in\mathbb{Z}$. That is, $f(x)$ is a null polynomial modulo
$p$. Since $\degree(f,p)<p$, one immediately derives that
$q_{e_1}(j_0)\equiv 0\pmod p$ for any integer
$j_0\in\{0,\cdots,p-1\}$. That is, $\forall e_1\in\{0,\cdots,p-1\}$,
$q_{e_1}(x)$ is a null polynomial modulo $p$.

When $n\geq 2$, assume this theorem is true for any integer less
than $n$. From Corollary \ref{corollary:base-null-polys-feature2},
$\widetilde{\mathcal{G}}_{p,i}(x)$ forms a function modulo $p$ of
period $p^{i+1}$. So, assuming $x^*=x\bmod
p^n=\sum_{k=n-1}^0j_kp^k$, one has
\[
\sum_{0\leq e_i\leq p-1 \atop 0\leq i\leq
n}\left(\prod_{i=1}^n\left(\widetilde{\mathcal{G}}_{p,i}(x)\right)^{e_i}q_{e_1,\cdots,e_n}(x)\right)\equiv\sum_{0\leq
e_i\leq p-1 \atop 0\leq i\leq
n}\left(\left(\widetilde{\mathcal{G}}_{p,n}(x)\right)^{e_n}\prod_{i=1}^{n-1}\left(\widetilde{\mathcal{G}}_{p,i}(x^*)\right)^{e_i}q_{e_1,\cdots,e_n}(x^*)\right)\equiv
0\pmod p.
\]
From Corollary \ref{corollary:base-null-polys-feature1}, $\forall
x^*\in\{0,\cdots,p^n-1\}$, if $i$ runs through a complete system of
residues modulo $p$, then $\widetilde{\mathcal{G}}_{p,n}(x)$ forms a
complete system of residues modulo $p$. This means that the
polynomial
\[
f(y)=\sum_{e_n=0}^{p-1}\left(\sum_{0\leq e_i\leq p-1 \atop 0\leq
i\leq
n-1}\left(\prod_{i=1}^{n-1}\left(\widetilde{\mathcal{G}}_{p,i}(x^*)\right)^{e_i}q_{e_1,\cdots,e_n}(x^*)\right)\right)y^{e_n}
\]
is congruent to 0 modulo $p$ for any integer $y\in\mathbb{Z}$. Then,
from $\degree(f,p)<p$, one immediately derives that
\[
\sum_{0\leq e_i\leq p-1 \atop 0\leq i\leq
n-1}\left(\prod_{i=1}^{n-1}\left(\widetilde{\mathcal{G}}_{p,i}(x^*)\right)^{e_i}q_{e_1,\cdots,e_n}(x^*)\right)\equiv
0\pmod p
\]
holds for any integer $x^*\in\{0,\cdots,p^n-1\}$ and
$e_n\in\{0,\cdots,p-1\}$. Now the value of $n$ is reduced to be
$n-1$, so one can use the hypothesis for each
$e_n\in\{0,\cdots,p-1\}$ to get the result: $\forall
e_1,\cdots,e_n\in\{0,\cdots,p-1\}$, $q_{e_1,\cdots,e_n}(x)$ is a
null polynomial modulo $p$. Thus, this theorem is proved.
\end{proof}

\begin{theorem}\label{theorem:H-is-least-degree-monic-null-poly}
Assume $p$ is a prime and $d\geq 1$. Then, $\mathcal{H}_{p,d}(x)$ is
a least-degree monic null polynomial modulo $p^d$.
\end{theorem}
\begin{proof}
Let us prove this theorem via induction on $d$. The case of $d=1$
has been proved above. Next, under the assumption that this lemma is
true for any integer less than $d$, let us prove the case of $d\geq
2$. Consider the following two different cases.

Case 1: $E_{\max}(\mathcal{H}_{p,d-1})<p$. Assume $f(x)$ is a null
polynomial modulo $p^d$. Following Corollary
\ref{corollary:null-poly-hierarchy-pd}, one has
$f(x)\equiv\sum_{j=d-1}^1p^{d-1-j}\mathcal{H}_{p,j}(x)q_j(x)+p^{d-1}q_0(x)\pmod{p^d}$,
where $q_j(x)$ is of degree less than $p$ except $q_{d-1}(x)$ and
$q_{j}(x)=0$ if $E_{\max}(\mathcal{H}_{p,j})=p$ (Lemma
\ref{lemma:monic-null-poly-degrees}). In other words, each exponent
of the effective polynomial $\mathcal{H}_{p,j}(x)$ satisfies $0\leq
e_{j,i}\leq p-1$. Then, from $f(x)\equiv 0\pmod{p^d}$ holds for any
integer $x$, dividing both sides by $p^{d-1}$, one has
$\sum_{j=d-1}^1\frac{\mathcal{H}_{p,j}(x)}{p^j}q_j(x)+q_0(x)\equiv
0\pmod p$, where $d_j$ is the sum of all exponents of
$\mathcal{H}_{p,j}(x)$. So we know that $\sum_{0\leq e_{j,k}\leq p-1
\atop 0\leq k\leq
n}\left(\prod_{k=1}^n\left(\widetilde{\mathcal{G}}_{p,k}(x)\right)^{e_{j,k}}q_{e_1,\cdots,e_n}(x)\right)\equiv
0\pmod p$ holds for any integer $x$, where $n$ is the maximum of the
digit lengths of $\langle\mathcal{H}_{p,j}\rangle$ for
$j\in\{1,\cdots,d-1\}$. Then, from Theorem
\ref{theorem:base-null-poly-hierarchy}, one immediately gets
$\forall e_{j,1},\cdots,e_{j,n}\in\{0,\cdots,p-1\}$,
$q_{e_1,\cdots,e_n}(x)$ is a null polynomial modulo $p$. That is,
$\forall j\in\{1,\cdots,p-1\}$, $q_j(x)$ is a null polynomial modulo
$p$. Then, choosing $p_{d-1}(x)=\NPp$, one can get a least-degree
monic null polynomial modulo $p^d$:
$f(x)=\mathcal{H}_{p,d-1}(x)\NPp=\mathcal{H}_{p,d}(x)$.

Case 2: $E_{\max}(\mathcal{H}_{p,d-1})=p$. From Lemma
\ref{lemma:monic-null-poly-degrees} and $\omega_1(p^{d+1})\geq
\omega_1(p^d)$, one can see $\mathcal{H}_{p,d}(x)$ must be a
least-degree monic null polynomial modulo $p^d$, since
$\degree(\mathcal{H}_{p,d},p)=\degree(\mathcal{H}_{p,d-1},p)$.

Combining the above two cases, this theorem is thus proved.
\end{proof}

\begin{lemma}
Assume $p$ is a prime and $d\geq 1$. Then,
$\omega_1(p^{d+1})-\omega_1(p^d)\in\{0,p\}$.
\end{lemma}
\begin{proof}
It is a direct result of the above theorem and Lemma
\ref{lemma:monic-null-poly-degrees}.
\end{proof}
\begin{corollary}
Assume $p$ is a prime and $d\geq 1$. Then, $p\mid \omega_1(p^d)$.
\end{corollary}
\begin{proof}
It is an obvious result of the above theorem.
\end{proof}

\begin{lemma}
Assume $p$ is a prime and $n\geq 1$, then
$\omega_1(p^{I_p(n)})=p^n$.
\end{lemma}
\begin{proof}
It is a direct result of Theorem
\ref{theorem:H-is-least-degree-monic-null-poly} and Lemma
\ref{lemma:G-is-null-polys}.
\end{proof}
\begin{theorem}\label{theorem:null-poly-omega1-pd-all}
Assume $p$ is a prime and $d\geq 1$. If
$\langle\mathcal{H}_{p,d}\rangle=(e_{d,n},\cdots,e_{d,1})_{I_p}$,
then
$\omega_1(p^d)=\sum_{i=1}^ne_{d,i}p^i=p\left(\sum_{i=1}^ne_{d,i}p^{i-1}\right)$.
\end{theorem}
\begin{proof}
It is a straightforward result of the above lemma and the definition
of $\mathcal{H}_{p,d}(x)$.
\end{proof}

\begin{algorithm}
From the above theorems and
$d=\langle\mathcal{H}_{p,d}\rangle=\sum_{i=n}^1e_{d,i}I_p(i)$, one
can derive an effective algorithm to determine all digits of
$\langle\mathcal{H}_{p,d}\rangle=(e_{d,n},\cdots,e_{d,1})_{I_p}$,
and then determine $\mathcal{H}_{p,d}(x)$ and the value of
$\omega_1(p^d)$ for any given $d\geq 1$. The algorithm can be
described in the following steps.
\begin{itemize}
\item
Step 1: find the integer $n\geq 0$ such that $I_p(n)\leq d<I_p(n+1)$
by calculating $n=\lceil\log_p(d(p-1)+1)\rceil-1$;

\item
Step 2: assign $d^{(n)}=d$, and then for $i=n\sim 2$, calculate
$e_{d,i}=\lfloor d^i/I_p(i)\rfloor$, $d^{(i-1)}=d^i-e_{d,i}I_p(i)$
and goto Step 3 if $d^{(i-1)}=0$;

\item
Step 3: assign $e_{d,1}=d^{(i-1)}$;

\item
Step 4: output
$\mathcal{H}_{p,d}(x)=\prod_{i=n}^1\left(\mathcal{G}_{p,i}(x)\right)^{e_{d,i}}$
and $\omega_1(p^d)=\sum_{i=1}^ne_{d,i}p^i$.
\end{itemize}
One can see that the time complexity and the space complexity of the
above algorithm are both
\[
O(n)=O\left(\lceil\log_p(d(p-1)+1)\rceil-1\right)=O(\log_pd).
\]
\end{algorithm}

\begin{example}
Assume $p$ is a prime and $i\geq 1$. Prove
$\mathcal{H}_{p,p^i}(x)=\left(\mathcal{G}_{p,i}(x)\right)^{p-1}\mathcal{G}_{p,1}(x)$
and $\omega_1(p^{p^i})=(p-1)p^i+p$.
\end{example}
\begin{solution}
From $I_p(n)=\frac{p^n-1}{p-1}$, one immediately gets $n=i$. Assign
$d^{(i)}=d$, one has $e_{d,i}=\lfloor
d^{(i)}/I_p(i)\rfloor=\left\lfloor\left.p^i\right/\left(\frac{p^i-1}{p-1}\right)\right\rfloor=p-1$
and $d^{(i-1)}=d^{(i)}-e_{d,i}I(i)=1$. Thus,
$e_{d,2}=\cdots=e_{d,i-1}=0$ and $e_{d,1}=1$. So
$\mathcal{H}_{p,p^i}(x)=\left(\mathcal{G}_{p,i}(x)\right)^{p-1}\mathcal{G}_{p,1}(x)$
and $\omega_1(p^{p^i})=(p-1)p^i+p$.
\end{solution}

Figure \ref{figure:general-list} gives a partial list of
$\mathcal{H}_{p,d}(x)$ from $d=1$ till $d=I_p(4)+p$. In the list,
the value of $\omega_1(p^d)$ increases by $p$ except at the end of
each row.
\begin{figure}[!htb]
\begin{equation*}
\begin{matrix}
& d=1\\
& \downarrow\\
\mathrm{N/A}^* & \mathcal{G}_{p,1}(x) & \cdots & (\mathcal{G}_{p,1}(x))^p\\
\hdashline \mathcal{H}_{p,I_p(2)}(x)=\mathcal{G}_{p,2}(x) &
\mathcal{G}_{p,2}(x)\mathcal{G}_{p,1}(x) & \cdots &
\mathcal{G}_{p,2}(x)(\mathcal{G}_{p,1}(x))^p\\
\vdots & \vdots & \ddots & \vdots\\
(\mathcal{G}_{p,2}(x))^{p-1} &
(\mathcal{G}_{p,2}(x))^{p-1}\mathcal{G}_{p,1}(x) & \cdots &
(\mathcal{G}_{p,2}(x))^{p-1}(\mathcal{G}_{p,1}(x))^p\\
(\mathcal{G}_{p,2}(x))^p\\
\hdashline \mathcal{H}_{p,I_p(3)}(x)=\mathcal{G}_{p,3}(x) &
\mathcal{G}_{p,3}(x)\mathcal{G}_{p,1}(x) & \cdots &
\mathcal{G}_{p,3}(x)(\mathcal{G}_{p,1}(x))^p\\
\mathcal{G}_{p,3}(x)\mathcal{G}_{p,2}(x) &
\mathcal{G}_{p,3}(x)\mathcal{G}_{p,2}(x)\mathcal{G}_{p,1}(x) &
\cdots &
\mathcal{G}_{p,3}(x)\mathcal{G}_{p,2}(x)(\mathcal{G}_{p,1}(x))^p\\
\vdots & \vdots & \ddots & \vdots\\
\mathcal{G}_{p,3}(x)(\mathcal{G}_{p,2}(x))^{p-1} &
\mathcal{G}_{p,3}(x)(\mathcal{G}_{p,2}(x))^{p-1}\mathcal{G}_{p,1}(x)
& \cdots &
\mathcal{G}_{p,3}(x)(\mathcal{G}_{p,2}(x))^{p-1}(\mathcal{G}_{p,1}(x))^p\\
\mathcal{G}_{p,3}(x)(\mathcal{G}_{p,2}(x))^p\\
(\mathcal{G}_{p,3}(x))^2 &
(\mathcal{G}_{p,3}(x))^2\mathcal{G}_{p,1}(x) & \cdots &
(\mathcal{G}_{p,3}(x))^2(\mathcal{G}_{p,1}(x))^p\\
(\mathcal{G}_{p,3}(x))^2\mathcal{G}_{p,2}(x) &
(\mathcal{G}_{p,3}(x))^2\mathcal{G}_{p,2}(x)\mathcal{G}_{p,1}(x) &
\cdots &
(\mathcal{G}_{p,3}(x))^2\mathcal{G}_{p,2}(x)(\mathcal{G}_{p,1}(x))^p\\
\vdots & \vdots & \ddots & \vdots\\
(\mathcal{G}_{p,3}(x))^{p-1}(\mathcal{G}_{p,2}(x))^{p-1} &
(\mathcal{G}_{p,3}(x))^{p-1}(\mathcal{G}_{p,2}(x))^{p-1}\mathcal{G}_{p,1}(x)
& \cdots &
(\mathcal{G}_{p,3}(x))^{p-1}(\mathcal{G}_{p,2}(x))^{p-1}(\mathcal{G}_{p,1}(x))^p\\
(\mathcal{G}_{p,3}(x))^{p-1}(\mathcal{G}_{p,2}(x))^p\\
(\mathcal{G}_{p,3}(x))^p\\
\hdashline \mathcal{H}_{p,I_p(4)}(x)=\mathcal{G}_{p,4}(x) &
\mathcal{G}_{p,4}(x)\mathcal{G}_{p,1}(x) & \cdots &
\mathcal{G}_{p,4}(x)(x)(\mathcal{G}_{p,1}(x))^p\\
\vdots & \vdots & \ddots & \vdots
\end{matrix}
\end{equation*}
\centering
\begin{minipage}{\textwidth}
$^*$ \footnotesize In fact, we can also generalize the result to
$d=0$. In this case, the modulus becomes $p^0=1$, so $\forall
x\in\mathbb{Z}$, one has $x\equiv 0\pmod 1$. Thus,
$f(x)=\prod_{i=1}^{\infty}(\mathcal{G}_{p,i})^0=1$ is a least-degree
monic null polynomial modulo 1. It is obvious that this case also
obeys the above theorems. In this sense, we can replace ``$d\geq 1$"
in the above theorems by ``$d\geq 0$".
\end{minipage}
\caption{An incomplete list of $\mathcal{H}_{p,d}(x)$, where $d$
increases by 1 from left to right and from top to bottom,
$\omega_1(p^d)$ does not change at the end of each row and increases
by $p$ elsewhere.}\label{figure:general-list}
\end{figure}

\subsection[Enumerating All Null Polynomials modulo $p^d$ ($d\geq 1$)]%
{Enumerating All Null Polynomials modulo $\bm{p^d}$ ($\bm{d\geq
1}$)} \label{subsection:null-poly-pd-enumerate-all}

Based on the fact that $\mathcal{H}_{p,d}(x)$ is a least-degree
monic null polynomial modulo $p^d$, one can enumerate all null
polynomials of (less than) a given degree modulo $p^d$.

\begin{theorem}\label{theorem:null-poly-enumerate-all}
Assume $p$ is a prime and $d\geq 1$. A polynomial $f(x)$ is a null
polynomial modulo $p^d$ if and only if
$f(x)\equiv\sum_{j=d}^1p^{d-j}\mathcal{H}_{p,j}(x)q_j(x)\pmod{p^d}$,
where
\begin{itemize}
\item
$q_d(x)$ is an arbitrary polynomial of arbitrary degree modulo $p$;

\item
when $j<d$, $q_j(x)=0$ if $E_{\max}(\mathcal{H}_{p,j})=p$;

\item
all other polynomials are arbitrary polynomials of degree less than
$p$ modulo $p$.
\end{itemize}
\end{theorem}
\begin{explain}
It can be easily derived from Theorem
\ref{theorem:H-is-least-degree-monic-null-poly}, following the
theorems on null polynomials modulo $p^d$ when $1\leq d\leq
p(p+1)+1$ (given in the previous subsection).
\end{explain}

The above theorem makes it possible to calculate the number of
(monic) null polynomials of a given degree modulo $p^d$, which may
be useful in some real applications. For example, through the
following theorem, the number of (monic) null polynomials is
actually the number of equivalent polynomials modulo $m$. This can
be further used to estimate the number of distinct polynomial
functions of some kind modulo $m$, once the number of candidate
polynomials of this kind have been known.

\begin{theorem}\label{theorem:Equ-Poly-Null-Poly}
Two polynomials, $f_1(x)$ and $f_2(x)$, are equivalent polynomials
modulo $m$ if and only if $f_1(x)-f_2(x)$ is a null polynomial
modulo $m$.
\end{theorem}
\begin{proof}
Obvious.
\end{proof}

\begin{definition}
In a complete system of polynomial residues modulo $p^d$, denote the
number of null polynomials of degree $n$ by $N_{np}(n,p^d)$ and the
number of monic null polynomials of degree $n$ by $N_{mnp}(n,p^d)$.
Similarly, denote the number of null polynomials of degree $\leq n$
by $N_{np}(\leq n,p^d)$ and the number of monic null polynomials of
degree $\leq n$ by $N_{mnp}(\leq n,p^d)$. Here, the subscript ``np"
means ``null polynomials" and ``mnp" denotes ``monic null
polynomials".
\end{definition}

It is obvious that $N_{np}(n,p^d)=N_{np}(\leq n,p^d)-N_{np}(\leq
n-1,p^d)$ and $N_{mnp}(n,p^d)=N_{mnp}(\leq n,p^d)-N_{mnp}(\leq
n-1,p^d)$. So, in the following we mainly focus on $N_{np}(\leq
n,p^d)$ and $N_{mnp}(\leq n,p^d)$.

\begin{theorem}
Assume $p$ is a prime, $d\geq 1$ and $n<\omega_1(p^d)$, then
$N_{mnp}(n,p^d)=N_{mnp}(\leq n,p^d)=0$.
\end{theorem}
\begin{proof}
It is a straightforward result of the definition of $\omega_1(p^d)$.
\end{proof}
\begin{theorem}
Assume $p$ is a prime and $d\geq 1$, then
$N_{mnp}(\omega_1(p^d),p^d)=N_{mnp}(\leq
\omega_1(p^d),p^d)=N_{np}(\leq\omega_1(p^d)-1,p^d)$.
\end{theorem}
\begin{proof}
The first equality is obvious since $N_{mnp}(\leq
\omega_1(p^d)-1,p^d)=0$. The second equality holds due to the
following fact: from Theorem \ref{theorem:null-poly-enumerate-all},
$q_d(x)\equiv 1\pmod{p^d}$ in a monic null polynomial of degree
$\omega_1(p^d)$, while $q_d(x)\equiv 0\pmod{p^d}$ in a null
polynomial of degree $\leq\omega_1(p^d)-1$ (all other coefficients
are free in both cases).
\end{proof}
\begin{theorem}
Assume $p$ is a prime, $d\geq 1$ and $n=\omega_1(p^d)+n^*$, where
$n^*\geq 0$. Then, it is true that
$N_{mnp}(n,p^d)=p^{dn^*}N_{mnp}(\omega_1(p^d),p^d)$ and
$N_{mnp}(\leq
n,p^d)=\frac{p^{d(n^*+1)}-1}{p^d-1}N_{mnp}(\omega_1(p^d),p^d)$.
\end{theorem}
\begin{proof}
From Theorem \ref{theorem:null-poly-enumerate-all}, to get a monic
null polynomial of degree $n$ modulo $p^d$, the highest coefficient
of $q_d(x)$ must be congruent to 1 modulo $p^d$ and all other
coefficients can be freely assigned. This means that
$N_{mnp}(n,p^d)=p^{dn^*}N_{mnp}(\omega_1(p^d),p^d)$. Then,
$N_{mnp}(\leq
n,p^d)=\sum_{i=0}^{n^*}p^iN_{mnp}(\omega_1(p^d),p^d)=\frac{p^{d(n^*+1)}-1}{p^d-1}N_{mnp}(\omega_1(p^d),p^d)$.
Thus this theorem is proved.
\end{proof}
\begin{theorem}
Assume $p$ is a prime, $d\geq 1$ and $n=\omega_1(p^d)+n^*$, where
$n^*\geq 0$, then $N_{np}(\leq
n,p^d)=p^dN_{mnp}(n,p^d)=\frac{p^{d(n^*+1)}(p^d-1)}{p^{d(n^*+1)}-1}N_{mnp}(\leq
n,p^d)$.
\end{theorem}
\begin{proof}
From Theorem \ref{theorem:null-poly-enumerate-all}, to get a monic
null polynomial $f(x)$ of degree $n$ modulo $p^d$, $q_d(x)$ should
be a monic polynomial of degree $n-\omega_1(p^d)$ modulo $p^d$; and
to get a null polynomial $f(x)$ modulo $p^d$, $q_d(x)$ can be
arbitrary polynomial of degree $\leq n-\omega_1(p^d)$ modulo $p^d$.
In other words, for the former case, the highest coefficient $a_n$
must be congruent to 1 modulo $p^d$, while for the latter case,
$a_n$ can be any value modulo $p^d$. Considering other $n-1$
coefficients can be freely assigned for both cases, one immediately
gets $N_{np}(\leq n,p^d)=p^dN_{mnp}(n,p^d)$. In a similar way and
from the above theorem, we can get $N_{np}(\leq
n,p^d)=\frac{p^{d(n^*+1)}}{(p^{d(n^*+1)})/(p^d-1)}N_{mnp}(\leq
n,p^d)=\frac{p^{d(n^*+1)}(p^d-1)}{p^{d(n^*+1)}-1}N_{mnp}(\leq
n,p^d)$. Thus this theorem is proved.
\end{proof}

\begin{remark}
In the above theorems, note that when $p^d$ is relatively large, we
have $N_{np}(\leq n,p^d)\approx p^dN_{mnp}(\leq n,p^d)$ and
$N_{mnp}(\leq n,p^d)\approx N_{mnp}(n,p^d)$. This means that
$N_{mnp}(n,p^d)\gg N_{mnp}(\leq n-1,p^d)$ in this case.
\end{remark}

\begin{theorem}
Assume $p$ is a prime, $d\geq 1$ and $n=\omega_1(p^d))+n^*$, where
$n^*\geq 0$. Then, $N_{np}(\leq
n,p^d)=p^{d(n^*+1)}N_{np}(\leq\omega_1(p^d)-1,p^d)$.
\end{theorem}
\begin{proof}
From Theorem \ref{theorem:null-poly-enumerate-all}, for each
polynomial counted in $N_{np}(\leq\omega_1(p^d)-1,p^d)$, there are
$p^{d(n^*+1)}$ possibilities of $q_d(x)$ modulo $p^d$.
\end{proof}

From the above theorem, one can only consider the value of
$N_{np}(\leq n,p^d)$ when $1\leq d\leq\omega_1(p^d)-1$.

\begin{theorem}\label{theorem:null-poly-numbers-d1d2}
Assume $p$ is a prime, $d_1,d_2\geq 1$ and
$n<\min(\omega_1(p^{d_1}),\omega_1(p^{d_2}))$. Then, $N_{np}(\leq
n,p^{d_1})=N_{np}(\leq n,p^{d_2})$.
\end{theorem}
\begin{proof}
It is a straightforward result of Theorem
\ref{theorem:null-poly-enumerate-all}, since $q_j(x)$ can be freely
assigned any value modulo $p^j$ for $j<
\min(\omega_1(p^{d_1}),\omega_1(p^{d_2}))$.
\end{proof}
\begin{corollary}
Assume $p$ is a prime, $d\geq 1$ and $n<\omega_1(p^d)$. If
$E_{\max}(\mathcal{H}_{p,d})=p$, then $N_{np}(\leq
n,p^d)=N_{np}(\leq n,p^{d+1})$.
\end{corollary}
\begin{proof}
It is result of the above theorem and the fact that
$\omega_1(p^d)=\omega_1(p^{d+1})$.
\end{proof}

\begin{definition}
Assume $p$ is a prime and $n\geq 1$. Define an integer sequence
$N_p(n)$ as follows:
\[
N_p(n)=\begin{cases}%
1, & \mbox{when }n=1,\\
\prod_{i=0}^{p-1}\left(p^{ip^{n-1}I_p(n-1)}N_p(n-1)\right)=p^{\frac{p^n\left(p^{n-1}-1\right)}{2}}(N(n-1))^p
& \mbox{when }n\geq 2.\\
\end{cases}
\]
Specially, $N_p(2)=p^{\frac{p^2(p-1)}{2}}$,
$N_p(3)=p^{\frac{p^3(p^2+p-2)}{2}}$ and
$N_p(4)=p^{\frac{p^4(p^3+p^2+p-3)}{2}}$.
\end{definition}

The above sequence $N(n)$ has an alternative (actually equivalent)
definition via addition (not product). In the following, we use both
of the two definitions to achieve a more concise description of
related results.

\begin{definition}
Assume $p$ is a prime and $n\geq 1$. Define an integer sequence
$\widetilde{N}_p(n)$ as follows:
\[
\widetilde{N}_p(n)=\begin{cases}%
0, & \mbox{when }n=1,\\
\sum_{i=0}^{p-1}\left(ip^{n-1}I_p(n-1)+\widetilde{N}_p(n-1)\right)
=\frac{p^n\left(p^{n-1}-1\right)}{2}+p\widetilde{N}_p(n-1), & \mbox{when }n\geq 2.\\
\end{cases}
\]
Then, define $N(n)_p=p^{\widetilde{N}_p(n)}$. Specially,
$\widetilde{N}_p(2)=\frac{p^2(p-1)}{2}$,
$\widetilde{N}_p(3)=\frac{p^3(p^2+p-2)}{2}$
$\widetilde{N}_p(4)=\frac{p^4(p^3+p^2+p-3)}{2}$.
\end{definition}

\begin{theorem}
Assume $p$ is a prime and $n\geq 1$, then
$\widetilde{N}_p(n)=\frac{p^n\left(\sum_{i=1}^{n-1}p^i-(n-1)\right)}{2}=\frac{p^n(I_p(n)-n)}{2}=\frac{p^n(p^n-np+(n-1))}{2(p-1)}$.
\end{theorem}
\begin{proof}
From the definition,
$\widetilde{N}_p(n)=\sum_{i=1}^{n-1}\left(p^{n-1-i}\left(\frac{p^{i+1}(p^i-1)}{2}\right)\right)=
\frac{p^n\left(\sum_{i=1}^{n-1}(p^i-1)\right)}{2}=\frac{p^n\left(\sum_{i=1}^{n-1}p^i-(n-1)\right)}{2}=
\frac{p^n}{2}\left(\sum_{i=0}^{n-1}p^i-n\right)=\frac{p^n(I_p(n)-n)}{2}=\frac{p^n}{2}\left(\frac{p^n-1}{p-1}-n\right)
=\frac{p^n(p^n-np+(n-1))}{2(p-1)}$. Thus this theorem is proved.
\end{proof}

\begin{lemma}
Assume $p$ is a prime and $n<p$, then $N_{np}(\leq n,p^d)=1$.
\end{lemma}
\begin{proof}
It is obvious, since $f(x)=0$ is the only null polynomial of degree
less than $p$.
\end{proof}
\begin{corollary}
Assume $p$ is a prime, then
$N_{np}(\leq\omega_1(p^{I_p(1)})-1,p^{I_p(1)})=1=N_p(1)$.
\end{corollary}
\begin{proof}
It is a straightforward result of the above lemma since
$\omega_1(p^{I_p(1)})-1=p-1<p$.
\end{proof}
\begin{lemma}
Assume $p$ is a prime, $1\leq d\leq p$ and $ip\leq n<(i+1)p$, where
$i\in\{1,\cdots,d-1\}$. Then $N_{np}(\leq
n,p^d)=p^{\frac{i(i-1)p}{2}+i(n-ip+1)}$. Specially,
$N_{np}(\leq\omega_1(p^d)-1,p^d)=p^{\frac{d(d-1)p}{2}}$ and
$N_{np}(\leq\omega_1(p^p)-1,p^p)=p^{\frac{p^2(p-1)}{2}}=N_p(2)$.
\end{lemma}
\begin{proof}
From Theorem \ref{theorem:null-poly-enumerate-all}, when $ip\leq
n<(i+1)p$, $q_1(x),\cdots,q_{i-1}(x)$ can be of degree $\leq p-1$
modulo $p^d$ but $q_i(x)$ should be of degree $\leq n-ip$ modulo
$p^d$. So, $N_{np}(\leq
n,p^d)=\sum_{j=1}^{i-1}p^{jp}+p^{i(n-ip+1)}=p^{\frac{i(i-1)p}{2}+i(n-ip+1)}$.
When $n=\omega_1(p^d)-1=pd-1$, one has $N_{np}(\leq
n,p^d)=p^{\frac{(d-1)(d-2)p}{2}+(d-1)((pd-1)-(d-1)p+1)}=p^{\frac{(d-1)(d-2)p}{2}+(d-1)p}=p^{\frac{d(d-1)p}{2}}$.
Further, when $d=p$, one immediately has
$N_{np}(\leq\omega_1(p^p)-1,p^p)=p^{\frac{p^2(p-1)}{2}}=N_p(2)$.
\end{proof}
\begin{corollary}
Assume $p$ is a prime, then
$N_{np}(\leq\omega_1(p^{I_p(2)})-1,p^{I_p(2)})=p^{\frac{p^2(p-1)}{2}}=N_p(2)$.
\end{corollary}
\begin{proof}
It is a result of the above lemma and the fact that
$\omega_1(p^{p+1})=\omega_1(p^p)$.
\end{proof}

\begin{theorem}
Assume $p$ is a prime and $n\geq 1$. Then,
$N_{np}(\leq\omega_1(p^{I_p(n)})-1,p^{I_p(n)})=N_p(n)$.
\end{theorem}
\begin{proof}
Let us prove this theorem via induction on $n$. The case of $n=1$ is
obvious and the case of $n=2$ has been proved above. Next, assume
this theorem is true for any integer less than $n$, let us prove the
case of $n\geq 3$.

From Theorem \ref{theorem:null-poly-enumerate-all},
$f(x)\equiv\sum_{j=d}^1p^{d-j}\mathcal{H}_{p,j}(x)q_j(x)\pmod{p^d}$.
When $d=I_p(n)$, we can rewrite this congruence as follows:
$f(x)\equiv\mathcal{G}_{p,n}(x)q_{I_p(n)}(x)+\sum_{0\leq e_{d,i}\leq
p-1 \atop 1\leq i\leq
n-1}p^{d-j(e_{d,1},\cdots,e_{d,n-1})}\prod_{i=1}^{n-1}\left(\mathcal{G}_{p,i}(x)\right)^{e_{d,i}}q_{j(e_{d,1},\cdots,e_{d,n-1})}(x)\pmod{p^d}$,
where $j(e_{d,1},\cdots,e_{d,n-1})=\sum_{i=1}^{n-1}e_{d,i}I(i)$. Let
us divide all the $q$-polynomials, excluding $q_{I_p(n)}(x)$, into
$p$ parts, each of which corresponds to a distinct value of
$e_{d,n-1}\in\{0,\cdots,p-1\}$:
\[
\left(\mathcal{G}_{p,n-1}(x)\right)^{e_{d,n-1}}\left(\sum_{0\leq
e_{d,i}\leq p-1 \atop 1\leq i\leq
n-2}p^{d-j(e_{d,1},\cdots,e_{d,n-1})}\prod_{i=1}^{n-2}\left(\mathcal{G}_{p,i}(x)\right)^{e_{d,i}}q_{j(e_{d,1},\cdots,e_{d,n-1})}(x)\right)\pmod{p^d}.
\]
One can see that the number of all possibilities of this part is
\begin{eqnarray*}
\prod_{0\leq e_{d,i}\leq p-1 \atop 1\leq i\leq
n-2}\left(p^{e_{d,n-1}I_p(n-1)}p^{j(e_{d,1},\cdots,e_{d,n-2})}\right)^p
& = & \left(p^{e_{d,n-1}pI_p(n-1)}\right)^{\omega_1\left(p^{I_p(n-2)}\right)}N_p(n-1)\\
& = & p^{e_{d,n-1}p^{n-1}I_p(n-1)}N_p(n-1),
\end{eqnarray*}
where note that $\omega_1(p^{I_p(n-2)})$ denotes the number of
factors in the left side. Finally, one immediately has
\begin{eqnarray*}
N_{np}(\leq\omega_1(p^{I_p(n)})-1,p^{I_p(n)}) & = &
\left(\prod_{e_{d,n-1}=0}^{p-1}\left(p^{e_{d,n-1}p^{n-1}I_p(n-1)}N_p(n-1)\right)\right)\\
& = & p^{\frac{p^n(p^{n-1}-1)}{2}}(N_p(n-1))^p=N_p(n).
\end{eqnarray*}
Thus this theorem is proved.
\end{proof}

\begin{definition}
Assume $p$ is a prime, $n\geq 1$ and $0\leq i\leq p$. Define a
two-index generalization of $N_p(n)$ as follows:
\[
N_p(n,i)=
\prod_{j=0}^{i-1}\left(p^{jp^nI_p(n)}N_p(n)\right)=p^{\frac{i(i-1)p^nI_p(n)}{2}}(N_p(n))^i.
\]
Specially, $N_p(n,0)=1$, $N_p(n,1)=N_p(n)$ and $N_p(n,p)=N_p(n+1)$.
\end{definition}

\begin{definition}
Similarly, define a two-index generalization of $\widetilde{N}_p(n)$
by
\[
\widetilde{N}_p(n,i)=\sum_{j=0}^{i-1}\left(jp^nI_p(n)+\widetilde{N}_p(n)\right)=\frac{i(i-1)p^nI_p(n)}{2}+i\widetilde{N}_p(n)=\log_p\left(N_p(n,i)\right).
\]
Specially, $\widetilde{N}_p(n,0)=0$,
$\widetilde{N}_p(n,1)=\widetilde{N}_p(n)$, and
$\widetilde{N}_p(n,p)=\widetilde{N}_p(n+1)$.
\end{definition}

\begin{theorem}
Assume $p$ is a prime, $n\geq 1$ and $d=iI(n)$, where $1\leq i\leq
p$. Then, $N_{np}(\leq\omega_1(p^d)-1,p^d)=N_p(n,i)$.
\end{theorem}
\begin{proof}
Following the same idea used in the proof of the above theorem, one
has
$N_{np}(\leq\omega_1(p^d)-1,p^d)=\sum_{e_{d,n}=0}^{i-1}\left(p^{e_{d,n}p^nI_p(n)}N_p(n)\right)
=p^{\frac{i(i-1)p^nI_p(n)(N_p(n))^i}{2}}=N_p(n,i)$.
\end{proof}
\begin{theorem}
Assume $p$ is a prime and $d\geq 1$. If
$\langle\mathcal{H}_{p,d}\rangle=(e_{d,n},\cdots,e_{d,1})_{I_p}$,
then
$N_{np}(\leq\omega_1(p^d)-1,p^d)=N_p(n,e_{d,n})\prod_{i=n-1}^1p^{d_i^*}N_p(i,e_{d,i})$,
where $d_i^*=\sum_{j=n}^{i+1}e_{d,j}p^jI_p(j)$.
\end{theorem}
\begin{proof}
This theorem can be easily proved following the same idea used in
the proofs of the above theorems. At first, we enumerate the number
of $q$-polynomials before
$\left(\mathcal{G}_{p,n}(x)\right)^{e_{d,n}}$, which is
$N_p(n,e_{d,n})$ from the above theorem. Then, we enumerate the
$q$-polynomials occurring between
$\left(\mathcal{G}_{p,n}(x)\right)^{e_{d,n}}$ and
$\left(\mathcal{G}_{p,n}(x)\right)^{e_{d,n}}\left(\mathcal{G}_{p,n-1}(x)\right)^{e_{d,n-1}}$,
which is $p^{e_{d,n}p^nI_p(n)}N_p(n-1,e_{d,n-1})$. Note that when
$e_{d,n-1}=0$ this number is still valid, though this case is not
covered by the above theorem. Repeat this procedure until all
$q$-polynomials except $q_d(x)$ are enumerated, we immediately prove
this theorem.
\end{proof}

\begin{remark}
When $n<\omega_1(p^d)-1$, the value of $N_{np}(\leq n,p^d)$ can be
calculated via an integer $d^*$ such that $\omega_1(d^*)-p\leq
n<\omega_1(d^*)$. Then, from Theorem
\ref{theorem:null-poly-numbers-d1d2}, $N_{np}(\leq
\omega_1(d^*)-1,p^{d^*})=N_{np}(\leq
n,p^{d^*})p^{\bar{d}^*(\omega_1(p^{d^*})-1-n)}=N_{np}(\leq
n,p^d)p^{\bar{d}^*(\omega_1(p^{d^*})-1-n)}$, where $\bar{d}^*$ is
the largest integer such that
$\omega_1(p^{\bar{d}^*})=\omega_1(p^{d^*})-1$. Thus, $N_{np}(\leq
n,p^d)=\frac{N_{np}(\leq
\omega_1(d^*)-1,p^{d^*})}{p^{\bar{d}^*(\omega_1(p^{d^*})-1-n)}}$.
\end{remark}

\section*{Appendix}

After finishing the first draft of this paper, we noticed that the
main result (Theorem
\ref{theorem:H-is-least-degree-monic-null-poly}) obtained in this
paper has been covered in Kempner's papers on this subject
\cite{Kempner:PolyResidue:TAMS1921a, Kempner:PolyResidue:TAMS1921b}
published in 1921. Kempner's proof is through a simple way that is
totally different from that one employed in this paper. In addition,
explicit formulas were not given to calculate $I_p(n)$ and
$\omega_1(p^d)$ (though a method is qualitatively explained in \S\S
1). Here, we give a brief introduction to Kempner's proof. Note that
we use definitions given in this paper to achieve a simpler
description of the results.

\begin{definition}
Define $\mu(m)$ to be the smallest positive integer such that
$\mu(m)!\equiv 0\pmod m$.
\end{definition}

\begin{lemma}[Lemma 1 in \cite{Kempner:PolyResidue:TAMS1921a}]
The polynomial $f(x)=\prod_{i=0}^{\mu(m)-1}(x-i)$ is a null
polynomial of degree $\mu(m)$ modulo $m$.
\end{lemma}
\begin{proof}
When $0\leq x\leq\mu(m)-1$, it is obvious $f(x)=0\equiv 0\pmod m$.
When $x\geq\mu(m)$, one has
$f(x)=\mu(m)!\binom{x}{\mu(m)}\equiv\mu(m)!\equiv 0\pmod m$, where
note that $\binom{x}{\mu(m)}$ is an integer since $x\geq\mu(m)$.
\end{proof}

\begin{lemma}[Lemma 2 in \cite{Kempner:PolyResidue:TAMS1921a}]
Any polynomial $f(x)=\sum_{k=0}^nc_kx^k$ is uniquely represented in
the form $\sum_{k=0}^na_k\binom{x}{k}$, and $c_n,a_n$ are both
different from zero if one of them is different from zero.
\end{lemma}
\begin{proof}
Since $x^k$ and $k!\binom{x}{k}$ are both monic polynomials of
degree $k$, we can divide $f(x)$ by $\{x^k\}_{k=n}^0$ to get the
first form and divide it by $k!\binom{x}{k}$ to get the second form.
\end{proof}
\begin{lemma}[Lemma 3 in \cite{Kempner:PolyResidue:TAMS1921a}]
In the above lemma, if $c_k$ ($k\in\{0,\cdots,n\}$) are integers,
then $a_k/k!$ and therefore $a_k$, are integers.
\end{lemma}
\begin{proof}
Comparing coefficients of the two forms of $f(x)$, this lemma is
immediately proved.
\end{proof}
\begin{lemma}[Lemma 4 in \cite{Kempner:PolyResidue:TAMS1921a}]
In the above lemma, if $f(x)$ is a null polynomial modulo $m$, then
$a_k\equiv 0\pmod m$ for $k\in\{0,\cdots,n\}$.
\end{lemma}
\begin{proof}
A straightforward result of the above lemma.
\end{proof}

\begin{lemma}[Lemma 5 in \cite{Kempner:PolyResidue:TAMS1921a}]
For a given modulus $m$, the integer polynomial
$f(x)=\prod_{i=0}^{\mu(m)-1}(x-i)$ is a monic least-degree null
polynomial modulo $m$, i.e., $\omega_1(m)=\mu(m)$.
\end{lemma}
\begin{proof}
Assume $f(x)$ is a null polynomial modulo $m$. Then, from the above
lemmas
\[
f(x)=\sum_{k=0}^nma_k\binom{x}{k}=ma_0+\sum_{k=1}^n\left(\frac{ma_k}{k!}\prod_{i=0}^{k-1}(x-i)\right),
\]
where $a_k$ are integers. To ensure $f(x)$ is a monic polynomial,
$\frac{ma_n}{n!}=1$, so $m\mid n!$. From the definition of $\mu(m)$,
we immediately have $n\geq\mu(m)$. Thus this lemma is proved.
\end{proof}

The above proof gives a different (slightly simpler than us)
least-degree monic null polynomial modulo $m$. In \S\S1 of
\cite{Kempner:PolyResidue:TAMS1921a}, Kempner discussed how to
calculate the value of $\mu(m)$ for different cases. When $m=p^d$
and $d\geq p$, he introduced an algorithm to calculate the value of
$\mu(m)$, which is actually the same as the one given in Algorithm 1
of this paper.

In addition, based on the above null polynomial, Kempner developed a
set of ``completely reduced polynomials" modulo $m$ and investigated
related problems about such completely reduced polynomials. For more
details, refer to Kempner's original work
\cite{Kempner:PolyResidue:TAMS1921a, Kempner:PolyResidue:TAMS1921b}.

\section*{Acknowledgments}

We would like to thank Prof. Carl Pomerance with the Department of
Mathematics, Dartmouth College, Prof. Melvyn B. Nathanson with the
Department of Mathematics, Lehman College, City University of New
York, and Prof. Zhijie Chen with the Department of Mathematics, East
China Normal University for some useful discussions. We also thank
Mr. Chengqing Li for sending us a copy of
\cite{Mullen:PolyFun-mod:AMH1984} for reference.

\bibliographystyle{unsrt}
\bibliography{Null-Poly}

\end{document}